\documentclass[11pt]{amsart}

\makeatletter
\usepackage{amssymb}
\usepackage{latexsym}
\usepackage{amsbsy}
\usepackage{amsfonts}

\def\marginpar#1{\ignorespaces}

\newtheorem{theorem}[equation]{Theorem}
\newtheorem{proposition}[equation]{Proposition}
\newtheorem{lemma}[equation]{Lemma}
\newtheorem{corollary}[equation]{Corollary}
\newtheorem{definition}[equation]{Definition}
\newtheorem{assumption}[equation]{Assumption}

\theoremstyle{definition}
\newtheorem{remark}[equation]{Remark}

\newtheorem{example}[equation]{Example}

\numberwithin{equation}{section}

\def\AArm{\fam0 \rm}%
\newdimen\AAdi%
\newbox\AAbo%
\def\AAk#1#2{\setbox\AAbo=\hbox{#2}\AAdi=\wd\AAbo\kern#1\AAdi{}}%

\newcommand{\BBone}{{\ensuremath{{\AArm 1\AAk{-.8}{I}I}}}}

\def\eqref#1{(\ref{#1})}
\def\eqlabel#1{\def\@currentlabel{#1}}

\def\formula#1{\def\@tempa{#1}\let\@tempb\theequation\def\theequation{%
\hbox{#1}}\def\@currentlabel{(\theequation)}$$}
\def\endformula{\leqno\hbox{(\@tempa)}$$\@ignoretrue\let\theequation\@tempb}

\def\given{\hskip5\p@\relax\vrule\@width.4\p@\hskip5\p@\relax}

\newcommand{\open}[1]{%
\par\normalfont\topsep6\p@\@plus6\p@\trivlist\item[\hskip\labelsep\itshape#1%
\@addpunct{.}]\ignorespaces}

\DeclareRobustCommand{\close}[1]{%
  \ifmmode % if math mode, assume display: omit penalty etc.
  \else \leavevmode\unskip\penalty9999 \hbox{}\nobreak\hfill
  \fi
  \quad\hbox{$#1$}}

\newlength{\toskip}

\def \R {{\mathbb R}}

\def \P {{\mathbb P}}
\def \E {{\mathbb E}}
\def \N {{\mathbb N}}

\def \L {{\mathbb L}}

\def \Var {\textrm{Var}}
\def \Ent {\textrm{Ent}}
\makeatother

\begin{document}
%%%%%%%%%%%%%%%%%%%%%%%%%%%%%%%%%%%%%%%%%%%%%%%%%%%%%%%%%%%%%%%%%%%%%%%
\date{\today}

\title[Lyapunov-Poincar\'e ...]{Rate of convergence for ergodic continuous Markov processes :
Lyapunov versus Poincar\'e.}

 \author[D. Bakry]{\textbf{\quad {Dominique} Bakry $^{\heartsuit}$ \, \, }}
\address{{\bf {Dominique} BAKRY},\\ Laboratoire de Statistiques et Probabilit\'es,
Universit\'e Paul Sabatier, 118 route de Narbonne  F- 31062 Toulouse cedex,  \\ and Institut
Universitaire de France} \email{bakry@math.ups-tlse.fr}

 \author[P. Cattiaux]{\textbf{\quad {Patrick} Cattiaux $^{\spadesuit}$ \, \, }}
\address{{\bf {Patrick} CATTIAUX},\\ Ecole Polytechnique, CMAP, F- 91128 Palaiseau cedex,
\\ and Universit\'e Paris X Nanterre, \'equipe MODAL'X, UFR SEGMI\\ 200 avenue de la
R\'epublique, F- 92001 Nanterre, Cedex.} \email{cattiaux@cmapx.polytechnique.fr}

 \author[A. Guillin]{\textbf{\quad {Arnaud} Guillin $^{\diamondsuit}$}}
\address{{\bf {Arnaud} GUILLIN},\\ Ecole Centrale Marseille et LATP \, Universit\'e  de Provence,
Technopole Château-Gombert, 39, rue F. Joliot Curie, 13453 Marseille Cedex 13.} \email{aguillin@egim-mrs.fr}

\maketitle
 \begin{center}
 \textsc{$^{\heartsuit}$ Universit\'e Paul Sabatier \quad and \quad Institut Universitaire de France}
\medskip

 \textsc{$^{\spadesuit}$ Ecole Polytechnique \quad and \quad Universit\'e Paris X}
\medskip

\textsc{$^{\diamondsuit}$ Ecole Centrale Marseille \quad and \quad Universite de Provence}
 \end{center}

\begin{abstract}
We study the relationship between two classical approaches for quantitative ergodic properties :
the first one based on Lyapunov type controls and popularized by Meyn and Tweedie, the second one
based on functional inequalities (of Poincar\'e type). We show that they can be linked through new
inequalities (Lyapunov-Poincar\'e inequalities). Explicit examples for diffusion processes are
studied, improving some results in the literature. The example of the kinetic Fokker-Planck
equation recently studied by H\'erau-Nier, Helffer-Nier and Villani is in particular discussed in
the final section.
\end{abstract}
\bigskip

\bigskip

\bigskip

\section{\bf Introduction, framework and first results.}\label{Intro}

Rate of convergence to equilibrium is one of the most studied problem in various areas of
Mathematics and Physics. In the present paper we shall consider a dynamics given by a time
continuous Markov process $(X_t, \P_x)$ admitting an (unique) ergodic invariant measure $\mu$, and
we will try to describe the nature and the rate of convergence to $\mu$.

In the sequel we denote by $L$ the infinitesimal generator (and $D(L)$ the extended domain of the generator),
by $P_t(x,.)$ the $\P_x$ law of $X_t$
and by $P_t$ (resp. $P_t^*$) the associated semi-group (resp. the adjoint or dual semi-group), so
that in particular for any regular enough density of probability $h$ w.r.t. $\mu$, $\int P_t(x,.)
h(x) \mu(dx) = P_t^* h \, d\mu$.
\medskip

Extending the famous Doeblin recurrence condition for Markov chains, Meyn and Tweedie developed
stability concepts for time continuous processes and furnished tractable methods to verify
stability \cite{MT2,MT3}. The most popular criterion certainly is the existence of a so called
Lyapunov function for the generator \cite{MT3,DMT}, yielding exponential (or geometric)
convergence, via control of excursions of well chosen functionals of the process. Sub-geometric or
polynomial convergence can also be studied (see \cite{Forro,Ver}
among others for the diffusion case). A very general form of the method is explained in the recent
work by Douc, Fort and Guillin \cite{DFG}, and we shall now explain part of their results in more
details.
\smallskip

\begin{definition}\label{defLyap}
Let $\phi$ be a positive function defined on $[1,+\infty[$. We say that $V \in D(L)$ is a
$\phi$-Lyapunov function if $V \geq 1$ and if there exist a constant $b$ and a closed petite set
$C$ such that for all $x$ $$LV(x) \, \leq \, - \, \phi(V(x)) \, + \, b \, \BBone_C(x) \, .$$

Recall that $C$ is a petite set if there exists some probability measure $a(dt)$ on $\R^+$ such
that for all $x \in C$ , $\int_0^{+\infty} P_t(x,.) a(dt) \geq \nu(.)$ where $\nu(.)$ is a
non-trivial positive measure.
\end{definition}

When $\phi$ is linear ($\phi(u)= a u$) we shall simply call $V$ a Lyapunov function. Existence of
a Lyapunov function furnishes exponential (geometric) decay \cite{DMT,MT3}, which is a particular
case of

\begin{theorem}\label{thmDFG}\cite{DMT} Thm 5.2.c, and \cite{DFG} Thm 3.10 and Thm 3.12.

Assume that there exists some increasing smooth and concave $\phi$-Lyapunov function $V$ such that
$V$ is bounded on the petite set $C$. Assume in addition that the process is irreducible in some
sense (see \cite{DMT,DFG} for precise statements). Then there exists a positive constant $c$ such
that for all $x$, $$\parallel P_t(x,.) - \mu \parallel_{TV} \, \leq \, c \, V(x) \, \psi(t) \, ,$$
where $\psi(t)=1/(\phi \circ H^{-1}_\phi) (t)$ for $H_\phi(t)= \int_1^t \, (1/\phi(s)) ds$, and
$\parallel . \parallel_{TV}$ is the total variation distance.
\end{theorem}

In particular if $\phi$ is linear, $\psi(t)=e^{-\rho t}$ for some positive explicit $\rho$.

Actually the result stated in Theorem \ref{thmDFG} can be reinforced by choosing suitable stronger
distances (stronger than the total variation distance actually weighted total variation distances)
but to the price of slower rates of convergence (see \cite{DMT,DFG} for details). In the same
spirit some result for some Wasserstein distance is obtained in \cite{HM}. An important drawback
of this approach is that there is no explicit control (in general) of $c$. One of the interest of
our approach will be to give explicit constants starting from the same drift condition.
\medskip

The pointwise Theorem \ref{thmDFG} of course extends to any initial measure $m$ such that $\int V
dm$ is finite. In particular, choosing $m=h \mu$ for some nice $h$, convergence reduces to the
study of $P_t^* h$ for large $t$. Long time behavior of Markov semi-groups is known to be linked
to functional inequalities. The most familiar framework certainly is the $\L^2$ framework and the
corresponding Poincar\'e (or weak Poincar\'e) inequalities, namely

\begin{theorem}\label{thmPoinc}
The following two statements are equivalent for some positive constant $C_P$
\begin{itemize}
\item[] \textbf{Exponential decay.} \quad  For all $f \in \L^2(\mu)$, $$\parallel P_t f - \int f
d\mu
\parallel_2^2 \, \, \leq \, e^{- t/C_P} \,  \parallel f - \int f d\mu \parallel_2^2 \, .$$
\item[] \textbf{Poincar\'e inequality.} \quad For all $f \in D_2(L)$ (the domain of the Fredholm
extension of $L$) , $$\Var_\mu(f) \, := \, \parallel f - \int f d\mu \parallel_2^2 \, \, \leq \,
C_P \, \int \, - 2f \, Lf \, d\mu \, .$$
\end{itemize}
In the sequel we shall define $\Gamma(f) = - 2f \, Lf$.
\end{theorem}

Thanks to Cauchy-Schwarz inequality and since $P_t$ and $P_t^*$ have the same $\L^2$ norm, a
Poincar\'e inequality implies an exponential rate of convergence in total variation distance, at
least for initial laws with a $\L^2$ density w.r.t. $\mu$.

As for the Meyn-Tweedie approach, one can get sufficient conditions for slower rates of
convergence, namely weak Poincar\'e inequalities introduced by Roeckner and Wang :

\begin{theorem}\label{thmWPI}\cite{r-w} Thm 2.1

Let $N$ be such that $N(\lambda f)=\lambda^2 N(f)$ , $N(P_t f) \leq N(f)$ for all $t$ and $N(f)
\geq
\parallel f
\parallel_2^2$.

Assume that there exists a non increasing function $\beta$ such that for all $s>0$ and all nice
$f$ the following inequality holds $$\textbf{(WPI)} \quad \parallel f - \int f d\mu
\parallel_2^2 \, \, \leq \, \beta(s) \, \int \Gamma(f) d\mu \, + \, s \, N\left(f - \int f
d\mu\right) \, .$$ Then
$$\parallel P_t f - \int f d\mu
\parallel_2^2 \, \, \leq \, \psi(t) \, N\left(f - \int f d\mu\right) \, ,$$ where $\psi(t)= 2 \inf \{s>0,
\beta(s) \log(1/s) \, \leq \,  t\}$.
\end{theorem}

In the symmetric case one can state a partial converse to Theorem \ref{thmWPI} (see \cite{r-w} Thm
2.3). Note that this time one has to assume that $N(h)$ is finite in order to get $\L^2$
convergence for $P_t^* h$. In general (WPI) are written with $N=\parallel . \parallel_\infty$ (or
the oscillation), criteria and explicit form of $\beta$ are discussed in \cite{r-w,BCR2}. A
particularly interesting fact is that any $\mu$ on $\R^d$ which is absolutely continuous w.r.t.
Lebesgue measure, $d\mu = e^{-F} dx$ with a locally bounded $F$, satisfies some (WPI).

Actually, as for the Meyn-Tweedie approach, one can show slower rates of convergence for less
integrable initial densities, as well as some results for an initial Dirac mass. We refer to
\cite{CGG} sections 4,5,6 for such a discussion in particular cases, we shall continue in this
paper. Actually \cite{CGG} is primarily concerned with (weak) logarithmic Sobolev inequalities,
that is replacing the $\L^2$ norm by the Orlicz norm associated to $u \rightarrow x^2 \log(1+x^2)$
i.e. replacing $\L^2$ initial densities by densities with finite relative entropy
(Kullback-Leibler information) w.r.t $\mu$. According to Pinsker-Csiszar inequality, relative
entropy dominates (up to a factor 2) the square of the variation distance, hence Gross logarithmic
Sobolev inequality (or its weak version introduced in \cite{CGG}) allows to study the decay to
equilibrium in total variation distance too.

Generalizations (interpolating between Poincar\'e and Gross) have been studied by several authors.
We refer to \cite{w00,w05,BCR1,BCR3,RZ} for related results on super-Poincar\'e and general
$F$-Sobolev inequalities, as well as their consequences for the decay of the semi-group in
appropriate Orlicz norms. We also refer to \cite{logsob} for an elementary introduction to the
standard Poincar\'e and Gross inequalities.
\medskip

If the existence of a $\phi$-Lyapunov function is a tractable sufficient condition for the
Meyn-Tweedie strategy (actually is a necessary and sufficient condition for the exponential case),
general tractable sufficient conditions for Poincar\'e or others functional inequalities are less
known (some of them will be recalled later), and in general no criterion is known (with the
notable exception of the one dimensional euclidean space). This is one additional reason to
understand the relationship between the Meyn-Tweedie approach and the functional inequality
approach, i.e. to link Lyapunov and Poincar\'e. This is the aim of this paper.
\medskip

Before to describe the contents of the paper, let us indicate another very attractive related
problem.

If $\mu$ is symmetric and ergodic, it is known that $\int \Gamma(f) d\mu = 0$ if and only if $f$
is a constant. In the non symmetric case this result is no more true, and we shall call fully
degenerate (corresponding to the p.d.e. situation) these cases.

Still in the symmetric case (or if $L$ is normal, i.e $L L^*=L^* L$), it is known that an
exponential decay $$\parallel P_t f - \int f d\mu
\parallel_2^2 \, \, \leq \, e^{- \rho t} \, N\left(f - \int f d\mu\right) \, ,$$ for some $N$ as
in Theorem \ref{thmWPI}, actually implies a (true) Poincar\'e inequality (see \cite{r-w} Thm 2.3).

A similar situation is no more true in the fully degenerate case. Indeed in recent works, H\'erau
and Nier \cite{HN04} and then Villani \cite{Vilhypo} have shown that for the kinetic Fokker-Planck
equation (which is fully degenerate) the previous decay holds with $N(g)=\parallel \nabla
g\parallel_2^2$ ($\mu$ being here a log concave measure, $N(g)$ is greater than the $\L^2$ squared
norm of $g$ up to a constant), and thanks to the hypoelliptic regularization property, it also
holds with $N(g)= C \, \Var_{\mu}(g)$ for some constant $C >1$ (recall that if $C\leq 1$, the
Poincar\'e inequality holds). Of course the Bakry-Emery curvature of this model is equal to
$-\infty$, otherwise an exponential decay with $N(g)= C \, \Var_{\mu}(g)$ would imply a Poincar\'e
inequality, even for $C>1$, so that this situation is particularly interesting.

It turns out that this model enters the framework of Meyn-Tweedie approach as shown in
\cite{wudamp} (also see \cite{DFG}). Hence relating Lyapunov to some Poincar\'e in such a case
(called hypocoercive by Villani) should help to understand the picture.

We shall also study this problem.
\bigskip

Let us briefly describe now our framework.

Recall that in all the paper $\mu$ is an invariant measure for the process with generator $L$.

The main additional hypothesis we shall make is the existence of a ``carr\'e du champ'', that is
we assume that there is an algebra which is a core for the generator and such that for $f$ and $g$
in this algebra
\begin{equation}\label{eqcarre}
L(fg) = f Lg + g Lf + \Gamma(f,g)
\end{equation}
where $\Gamma(f,g)$ is the polarization of $\Gamma(f)$. We shall also assume that $\Gamma$ comes
from a derivation, i.e. for $f$, $h$ and $g$ as before
\begin{equation}\label{eqderive}
\Gamma(fg,h) \, = \, f \Gamma(g,h) + g \Gamma(f,h) \, .
\end{equation}
The meaning of these assumptions in terms of the underlying stochastic process is explained in the
introduction of \cite{cat4}, to which the reader is referred for more details (also see
\cite{bakry} for the corresponding analytic considerations).

Applying Ito's formula, we then get that for all smooth $\Psi$, and $f$ as before,
\begin{equation}\label{eqchainrule}
L \Psi(f) \, = \, \frac{\partial \Psi}{\partial x}(f) \, Lf  \, + \, 1/2 \, \frac{\partial^2
\Psi}{\partial x^2}(f) \, \Gamma(f) \, .
\end{equation}

\bigskip
Our plan will be the following: in the second Section we show how to get controls in variance or
in entropy starting from the result of Theorem 1.2 which will be seen to be quite sharp. The
Section 3 will be devoted to the introduction of (weak) Lyapunov Poincar\'e inequalities, leading
to tractable criteria enabling us to give explicit control of convergence via ($\phi$-)Lyapunov
condition, illustrated by the examples of Section 4. The next section presents similar results for
the entropy, before presenting in the final Section an application in the particular fully
degenerate case.

%%%%%%%%%%%%%%%%%%%%%%%%%%
%%%%%%%%%%%%%%
%%%%%%%%%%%%%%%%%%%%%%%%%%

\section{\bf From Lyapunov to Poincar\'e and vice versa.}\label{MeynTweedie}

We first show that, in the symmetric case, the Meyn-Tweedie method immediately furnishes some
Poincar\'e inequalities.

Indeed let us assume that the hypothesis of Theorem \ref{thmDFG} are fulfilled, and let $f$ be a
bounded function such that $\int f d\mu = 0$. Then, if $f$ does not vanish identically, we may
define $h=f_+/\int f_+ d\mu$ which is a bounded density of probability. Thus, if $V \in
\L^1(\mu)$, $\int h V d\mu < +\infty$. It follows that $\parallel P_t^*h - 1
\parallel_{\L^1(\mu)}$, which is the total variation distance between $\mu$ and the law at time
$t$ (starting from $h \mu$), goes to 0 as $t \to +\infty$, with rate $c \psi(t)$ defined in
Theorem \ref{thmDFG}.

Hence for $0<\beta < 1$,
\begin{eqnarray*}
\int |P_t^*f_+ - \int f_+ d\mu|^2 \, d\mu & = & \left(\int f_+ d\mu\right)^2 \,
\int \, \left(P_t^*h - 1\right)^2 \, d\mu \\
& \leq & \left(\int f_+ d\mu\right)^2 \, \int \, \left(P_t^*h - 1\right)^\beta \,  \left(P_t^*h -
1\right)^{2-\beta} \, d\mu \\ & \leq & \left(\int f_+ d\mu\right)^2 \, \left(\int |P_t^*h - 1|
d\mu\right)^{\beta} \, \left(\int |P_t^*h - 1|^{\frac{2-\beta}{1-\beta}} d\mu\right)^{1-\beta} \\
& \leq & c_\beta \, \psi^{\beta}(t) \, \left(\int f_+ V d\mu\right)^\beta \, \left(\int |f_+ -
\int f_+ d\mu|^{\frac{2-\beta}{1-\beta}} d\mu\right)^{1-\beta} \\ & \leq & c_\beta \,
\psi^{\beta}(t) \, \left(\int f_+ V d\mu\right)^\beta \, \left(\int
(2f_+)^{\frac{2-\beta}{1-\beta}} d\mu\right)^{1-\beta}
\end{eqnarray*}
where we have used that $P_t^*$ is an operator with norm equal to 1 in all the $\L^p$'s ($p\geq
1$), the elementary $|a+b|^p \leq 2^{p-1}(|a|^p+|b|^p)$ for $p\geq 1$ and H\"{o}lder inequality.

Thus

\begin{eqnarray*}
\int (P_t^* f)^2 d\mu & \leq &  2 \int |P_t^*f_+ - \int f_+ d\mu|^2 \, d\mu \, + \, 2 \int
|P_t^*f_- - \int f_- d\mu|^2 \, d\mu \\ & \leq & 2^{4-\beta} \, c_\beta \, \psi^{\beta}(t) \,
\left(\int |f| V d\mu\right)^\beta \, \left(\int |f|^{\frac{2-\beta}{1-\beta}}
d\mu\right)^{1-\beta} \, .
\end{eqnarray*}

In the symmetric case (or more generally the normal case) we may thus apply Theorem 2.3 in
\cite{r-w} so that we have shown

\begin{theorem}\label{thmMT}
Under the hypotheses of Theorem \ref{thmDFG}, for any $f$ such that $\int f d\mu =0$ and any
$0<\beta<1$ it holds $$\int (P_t^* f)^2 d\mu \, \leq \, C_\beta \, \psi^{\beta}(t) \, \left(\int
|f| V d\mu\right)^\beta \, \left(\int |f|^{\frac{2-\beta}{1-\beta}} d\mu\right)^{1-\beta} \, .$$
The result extends to $\beta=1$ provided $f$ is bounded, and in this case $$\int |P_t^* f|^2 d\mu
\, \leq \, C \, (\int V d\mu) \, \parallel f\parallel_\infty^2 \, \psi(t) \, .$$

If in addition $\mu$ is a symmetric measure for the process, then $\mu$ satisfies a Weak
Poincar\'e Inequality with $N(f)= C(V) \, \parallel f\parallel_\infty^2$ and $$\beta(s) = s \,
\inf_{u>0}\, {1\over u}\,\psi^{-1}(u \, e^{(1-u/s)}) \, \textrm{ with } \, \psi^{-1}(a) := \inf
\{b>0 \, , \, \psi(b)\leq a\} \, .$$ In particular if $\psi(t)=e^{-\rho t}$, $\mu$ satisfies a
Poincar\'e inequality.
\end{theorem}

The fact that a Lyapunov condition furnishes some Poincar\'e inequality in the symmetric case is
already known see Wu \cite{wu1,wu2}, but the techniques used by Wu are different and rely mainly
on spectral ideas. Note also that a Lyapunov function is always in $\L^1(\mu)$ by integrating the
Lyapunov condition, otherwise only $\phi\circ V$ is integrable w.r.t. $\mu$ (as a direct
consequence of the Lyapunov inequality). In fact, the simple use, of this theorem enables us to
derive very easily the correct rate of convergence to equilibrium and to extend known sharp weak
Poincar\'e inequality in dimension one to higher dimension. The major drawback is that the
constants are quite unknown in the general case, however we refer to \cite{DMR} to results
providing explicit constants.

The same idea furnishes without any effort a similar result for the decay of relative entropy.
Indeed, if $h$ is a density of probability (w.r.t. $\mu$), using the concavity of the logarithm,
we get for any $0<\beta<1$
\begin{eqnarray*}
\int P^*_t h \, \log P^*_t h \, d\mu & = & \int (P^*_t h \, - \, 1) \, \log P^*_t h \, d\mu \, +
\, \int \log P^*_t h \, d\mu \\ & \leq & \int (P^*_t h \, - \, 1) \, \log P^*_t h \, d\mu \, + \, \log
\left( \int  P^*_t h \, d\mu\right) \\ & \leq & \int |P^*_t h \, - \, 1| \, |\log P^*_t h| \, d\mu
\\ & \leq & \left(\int |P^*_t h \, - \, 1| d\mu\right)^{\beta} \, \left(\int |P^*_t h \, - \, 1| \,
|\log P^*_t h|^{\frac{1}{1-\beta}} \, d\mu \right)^{1-\beta} \, .
\end{eqnarray*}
It is easily seen that the function $u \mapsto |u-1| |\log u|^p$ is convex on $]0,+\infty[$ for
$p\geq 1$, so that $$t \, \mapsto \, \int |P^*_t h \, - \, 1| \, |\log P^*_t
h|^{\frac{1}{1-\beta}} \, d\mu$$ is decaying on $\mathbb R^+$. We thus have obtained

\begin{theorem}\label{thmMTlog}
Under the hypotheses of Theorem \ref{thmDFG}, for any non-negative $h$ such that $\int h d\mu =1$
and any $0<\beta<1$ it holds $$\int P^*_t h \, \log P^*_t h \, d\mu \, \leq \, C_\beta \,
\psi^{\beta}(t) \, \left(\int h V d\mu\right)^\beta \, \left(\int |h-1| |\log
h|^{\frac{1}{1-\beta}} d\mu\right)^{1-\beta} \, .$$ The result extends to $\beta=1$ provided $h$
is bounded, and in this case $$\int P^*_t h \, \log P^*_t h \, d\mu \, \leq \, C \, (\int V d\mu)
\, \parallel h\parallel_\infty \,  \log (\parallel h\parallel_\infty) \, \psi(t) \, .$$
\end{theorem}

Note that (in the symmetric case) there is no analogue converse result for relative entropy as for
the variance. Indeed recall that if $h$ is a density of probability, $\int h \log h d\mu \leq
\Var_\mu(h)$, hence relative entropy is decaying exponentially fast, controlled by the initial
variance of $h$ as soon as a Poincar\'e inequality holds. But it is known that a Poincar\'e
inequality may hold without log-Sobolev inequality. However, starting from Theorem \ref{thmMTlog}
one can prove some (loose) weak log-Sobolev inequality, see \cite{CGG} sections 4 and 5.
\medskip

Of course Theorem \ref{thmMT} and Theorem \ref{thmMTlog} furnish (in the non-symmetric case as
well) controls depending on the integrability of $V$. For instance if $V$ has all polynomial
moments, we may control $\int |f| V d\mu$ by some $\int |f|^p d\mu$ in Theorem \ref{thmMT} and if
$\int e^{qV} d\mu < +\infty$ for some $q>0$ we may control $\int h V d\mu$ by the $u \log_+ u$
Orlicz norm of $h$ in Theorem \ref{thmMTlog}. Recall that a Lyapunov function is in $\L^1(\mu)$.
\medskip

We have seen that, in the symmetric case, the existence of a Lyapunov function implies a
Poincar\'e inequality. Let us briefly discuss some possible converse.

If $P_t$ is $\mu$ symmetric for some $\mu$ satisfying a Poincar\'e inequality, then we know that
$P_t$ has a spectral gap, say $\theta$. Let $f$ be an eigenfunction associated with the eigenvalue
$-\theta$, i.e. $L f + \theta f = 0$. If the semi-group is regularizing (in the ultracontractive
case for instance), $f$ has to be bounded. Assume that $f$ is actually bounded and say continuous.
Since $\int f d\mu=0$, changing $f$ into $-f$ if necessary, we may assume that $\sup f \geq - \inf
f= - M$. Then define $g=f+1+M$. $Lg = - \theta g + \theta (1+M)$, so that for all $0<\kappa<1$,
$$Lg \, \leq \, - \kappa \, \theta \, g \, + \, \theta \, \kappa \, (1+M) \, \BBone_C$$ with $C=
\{f \leq (1+M)\kappa/(1-\kappa)\}$ a non empty (and non full) closed set.

Of course the previous discussion only covers very few cases, but it indicates that some converse
has to be studied.

\medskip

Another possible way to prove a converse result is the following. Assume that
$d\mu(x)=e^{-2V(x)}dx$ where $V$  is $C^3$ and such that \begin{equation}\label{conditionV}
|\nabla V|^2(x) \, - \, \Delta V(x) \, \geq \, - C_{min} \, > - \infty
\end{equation}
 for a nonnegative
$C_{min}$ so that the process defined by (recalling that $B_t$ is an usual Brownian motion in $\R^d$)
\begin{equation}
\label{eds}
dX_t \, = \, dB_t \, -  \, (\nabla V)(X_t)dt \quad , \quad Law(X_0)=\nu
\end{equation}
has a unique non explosive strong solution. Assume also that $\mu$ satisfies a Poincar\'e
inequality.  The difficulty here is that by using Poincar\'e inequality we inherit a control for
all smooth $f$ with finite variance as $$\Var_\mu(P_tf)\le e^{-\lambda t}\Var_\mu(f).$$ But a
drift inequality concerns the generator and its behavior towards some chosen function for all $x$.
However it is known, see Down-Meyn-Tweedie \cite{DMT} (Th. 5.2, Th. 5.3 and the remarks after Th.
5.3), that the existence of a drift condition is ensured by $$\|P_t\delta_x-\pi\|_{TV}\le
M(x)\rho^t$$ for some larger than 1 function $M$ and $\rho<1$. But it is once again a control
local in $x$. In this direction, one can show (see \cite[Theorem 3.2.7]{Ro99}) that
$\Ent_{\mu}{P_t \delta_x}$ is finite for all $t>0$. But control in entropy is not useful as our
assumption is a Poincar\'e inequality and thus a control in $L^2$ is needed. Actually the proof of
Royer can be used in order to get the following result. Replacing the convex $\gamma$ therein by
$\gamma(y)=y^2$ we obtain $$\int (P_t \delta_x)^2 d\mu \leq Z \, e^{2V(x)} \, \E\left[e^{-2v(B_t)}
\, e^{- \frac 12 \, \int_0^t [|\nabla V|^2 - \Delta V](B_s) ds} \right] \leq Z \, e^{2V(x)} \,
e^{\frac 12 \, C_m t} $$ where $e^{-2v(y)} = (2\pi t)^{-d/2} \, e^{- |y-x|^2/2t}$. By the
Poincar\'e inequality, we then get that for some $\lambda$ and $t_0$
$$\Var_\mu(P_t\delta_x)\le e^{-\lambda (t-t_0)}\Var_\mu(P_{t_0}\delta_x)\le Z
\, e^{2V(x)} \, e^{\frac 12 \, C_m t_0}\,e^{-\lambda (t-t_0)}$$ which ends the work as a control
in $L^2$ enables us to control the $L^1$ distance, and we thus get the existence  of a Lyapunov
function. However, the Lyapunov function $V$ is not available in close form (see \cite{MT,DMT} for
a precise formula).

\medskip

Finally, let us mention that it is not possible to get a converse result as previously starting
from a weak  Poincar\'e inequality as 1) we do not know how to control $\|P_t\delta_x\|_\infty$
(even if it should be controlled in many case) and 2) there is no converse part in the
Meyn-Tweedie framework (even in the discrete time case) for sub geometric convergence in total
variation towards $\phi$-Lyapunov condition.

\bigskip

%%%%%%%%%%%%%%%%%%%%%%%%%%%%%%%%%%%%%%%%%%%%%%%%%%%%%%%%%%%%%%%%%%%%%%%%%
%%%%
%%%%section 3
%%%%
%%%%%%%%%%%%%%%%%%%%%%%%%%%%%%%%%%%%%%%%%%%%%%%%%%%%%%%%%%%%%%%%%%%%%%%%%

\section{\bf From Lyapunov to Poincar\'e. Continuation.}\label{LyapunovPoincare}

Since we have seen in the previous section that the existence of Lyapunov functions furnishes
functional inequalities in the symmetric case, in this section we shall study relationship between
some modified Poincar\'e inequality (still yielding exponential decay) and the existence of a
Lyapunov function (with $\phi(u)=\alpha u$), without assuming symmetry.

\subsection{\bf Lyapunov-Poincar\'e inequalities.}\label{subsLyapunov}

The key tool is the following elementary lemma

\begin{lemma}\label{lemmecle}
For $\Psi$ smooth enough, $W \in D(L)$ and $f \in L^{\infty}$, define $I^\Psi_W(t)= \int \,
\Psi(P_t f) \, W \, d\mu$. Then for all $t>0$,
$$\frac{d}{dt} \, I^\Psi_W(t) \, = \, - \, \int \, 1/2 \, \Psi''(P_t f) \, \Gamma(P_t f) \, W \, d\mu \,
+ \, \int \, L^*W \, \Psi(P_t f) \, d\mu \, .$$ In particular for $\Psi(u)=u^2$ we get (denoting
simply by $I_W$ the corresponding $I_W^\Psi$)
$$I'_W(t) \, = \, - \, \int \, \Gamma(P_t f) \, W \,
d\mu \, + \, \int \, L^*W \, P_t^2 f \, d\mu \, .$$
\end{lemma}
\begin{proof}
Recall that $\int L(\Psi(g) W) d\mu =0$. Using \eqref{eqcarre} and \eqref{eqchainrule} with $g=P_t
f$ we thus get
\begin{eqnarray*}
\frac{d}{dt} \, I^\Psi_W(t) & = & \int \, \Psi'(P_t f) \, LP_t f \, W \, d\mu \, \\ & = & \int
\left( L (\Psi(P_t f)) - 1/2 \Psi''(P_t f) \Gamma(P_t f) \right) \, W \, d\mu
\end{eqnarray*}
hence the result.
\end{proof}

This Lemma naturally leads to the following Definition and Proposition
\begin{definition}\label{defLP}
We shall say that $\mu$ satisfies a (W)-Lyapunov-Poincar\'e inequality, if there exists $W\in
D(L)$ with $W\geq 1$ and a constant $C_{LP}$ such that for all nice $f$ with $\int f d\mu=0$,
$$\int \, f^2 \, W \, d\mu \, \leq \, C_{LP} \, \int \, \left(W \,
\Gamma(f) \, - \, f^2 \, LW\right) \, d\mu \, .$$
\end{definition}

%\hrule
%{\bf Comments from AG}
%this inequality is not stable by pertubation (due to the mean...).... that can be a problem for application...???
% unless we explain that we do it for nice measure to get weak poincare which is stable...
%\smallskip
%\hrule

\begin{proposition}\label{proptrou}
The following statements are equivalent \begin{itemize} \item $\mu$ satisfies a (W)-
Lyapunov-Poincar\'e inequality, \item $\int (P_t^*f)^2  \, W \, d\mu \, \leq \, e^{- (t/C_{LP})}
\, \int f^2 W d\mu$ for all $f$ with $\int f d\mu = 0$. \end{itemize} In particular for all $f$
such that $\int f^2 W d\mu < +\infty$, $P_t f$ and $P_t^* f$ go to $\int f d\mu$ in $\L^2(\mu)$
with an exponential rate.
\end{proposition}
\begin{proof}
We consider $I^*_W(t)$ replacing $P_t$ by $P^*_t$. Taking the derivative at time $t=0$ furnishes
as usual the converse part. For the direct one, we only have to use Gronwall's lemma. Indeed the
Lyapunov-Poincar\'e inequality yields $(I^*_W)'(t) \, \leq - \, (1/C_{LP}) \, I^*_W(t)$. Since
$I^*_W(t)$ is non-negative, this shows that $I^*_W(t)$ is non increasing, hence converges to some
limit as $t$ tends to infinity, and this limit has to be 0 (otherwise $I^*_W$ would become
negative). Since $I^*_W(+\infty)=0$, the result follows by integrating the differential inequality
above.
\end{proof}

Note that a Lyapunov-Poincar\'e inequality is not a weighted Poincar\'e inequality (we still
assume that $\int f d\mu=0$) and depends on the generator $L$ (not only on the carr\'e du champ).
But as we already mentioned, Theorem 2.3 in \cite{r-w} tells us that, in the symmetric case, if
$$\int P_t^2 f \, d\mu \, \leq \, c \, e^{-\delta t} \, \parallel f \parallel_\infty^2$$ for all $f$
such that $\int f d\mu =0$, then $\mu$ satisfies the usual Poincar\'e inequality, with
$C_P=1/\delta$. Hence
\begin{corollary}\label{corpoinc}
If $L$ is $\mu$ symmetric and $\mu$ satisfies a (W)-Lyapunov-Poincar\'e inequality for some $W \in
\L^1$, then $\mu$ satisfies the ordinary Poincar\'e inequality, with $C_P=C_{LP}$.
\end{corollary}

Now we turn to sufficient conditions for a Lyapunov-Poincar\'e inequality to hold.

Recall that we called $V$ a Lyapunov function if $LV \, \leq - \alpha V + b \BBone_C$. Note that
integrating this relation w.r.t. $\mu$ yields $\alpha \int V d\mu \leq b \mu(C)$, so that, first
we have to assume that $\int V d\mu < +\infty$, second since $V\geq 1$, $b$ and $\mu(C)$ have to
be positive.

Before stating the first result of this section we shall introduce some definition

\begin{definition}\label{defpoincloc}
Let $U$ be a subset of the state space $E$. We shall say that $\mu$ satisfies a local Poincar\'e
inequality on $U$ if there exists some constant $\kappa_U$ such that for all nice $f$ with $\int_E
\, f d\mu =0$,
$$\int_U \, f^2 \, d\mu \, \leq \, \kappa_U \, \int_E \Gamma(f) d\mu \, + \, (1/\mu(U)) \,
\left(\int_U \, f \, d\mu\right)^2 \, .$$
\end{definition}

Notice that the energy integral in the right hand side is taken over the whole space $E$. We may
now state

\begin{theorem}\label{thmlyapgap}
Assume that there exists a Lyapunov function $V$ i.e. $LV \leq - 2 \alpha V + b \BBone_C$ for some
set $C$ (non necessarily petite).

Assume that one can find a (large) set $U$ such that $\mu$ satisfies a local Poincar\'e inequality
on $U$.

Assume in addition that
\begin{enumerate}
\item[(1)] either $U$ contains $C'=C \cap \{V \leq b/\alpha\}$ and $\alpha \, \mu(U) > b \mu(U^c)$
, \item[(2)] or $U$ contains $\{V \leq b/\alpha\}$ and $\mu(U)>\mu(U^c)$.
\end{enumerate}
Then one can find some $\lambda>0$ such that if $W=V+\lambda$, $\mu$ satisfies a
(W)-Lyapunov-Poincar\'e inequality.

More precisely, corresponding to the two previous cases one can choose
\begin{enumerate}
\item[(1)] $\lambda = (b \kappa_U - 1)_+$ and $1/C_{LP}=\alpha \left(1 - \frac{b \,
\mu(U^c)}{\alpha \, \mu(U)}\right)/(1+\lambda)$, \item[(2)] or $\lambda = (b \kappa_U - 1)_+$ and
$1/C_{LP}=\alpha \left(1 - \frac{\mu(U^c)}{\mu(U)}\right)/(1+\lambda)$.
\end{enumerate}
\end{theorem}

\begin{proof}
First remark the following elementary fact : define $C'=C \cap \{V \leq b/\alpha\}$. Then $LV \leq
- \alpha V + b \BBone_{C'}$, that is we can always assume that $C$ is included into some level set
of $V$. In the sequel $\theta=b/\alpha$. First we assume that $U$ contains $\{V \leq b/\alpha\}$,
so that it contains $C'$.

Let $\int f d\mu = 0$. Then for all $\lambda >0$ it holds
\begin{eqnarray*}
\int f^2 \, (V+\lambda)\, d\mu & \leq & (1+\lambda) \, \int f^2 \, V \, d\mu \\ & \leq &
(1+\lambda)/\alpha \, \int f^2 \, \left(-L(V+\lambda) + b \BBone_{C'}\right) \, d\mu \, .
\end{eqnarray*}
But since $\int_U f d\mu = - \int_{U^c} f d\mu$ it holds
\begin{eqnarray*}
\int_{C'}  \, f^2 \, d\mu & \leq  & \int_U \, f^2 \, d\mu \, \leq \, \kappa_U \, \int \Gamma(f)
d\mu \, + \, (1/\mu(U)) \, \left(\int_U \, f \, d\mu\right)^2 \\ & \leq & \kappa_U \, \int
\Gamma(f) d\mu \, + \, (1/\mu(U)) \, \left(\int_{U^c} \, f \, d\mu\right)^2 \\ & \leq & \kappa_U
\, \int \Gamma(f) d\mu \, + \, (\mu(U^c)/\mu(U)) \, \left(\int_{U^c} \, f^2 \, d\mu\right) \\ &
\leq & \kappa_U \, \int \Gamma(f) d\mu \, + \, (\mu(U^c)/\theta \, \mu(U)) \, \left(\int_{U^c} \,
f^2 \, V \, d\mu\right)\, ,
\end{eqnarray*}
where we used $V/\theta \geq 1$ on $U^c$. So, if we choose $\lambda=(b \kappa_U - 1)_+$ we get
$$b \, \int_{C'} \, f^2 \, d\mu \, \leq \, \int \, \Gamma(f) \, (V+\lambda) \, d\mu \, + \, (b
\mu(U^c)/\theta \, \mu(U)) \, \left(\int \, f^2 \, V \, d\mu\right) \, .$$ It yields $$\int \,
\left(W \, \Gamma(f) \, - \, f^2 \, LW\right) \, d\mu \, \geq \, \alpha \, \left(1 - \frac{b \,
\mu(U^c)}{\theta \, \alpha \, \mu(U)}\right) \, \int f^2 \, V \, d\mu \, ,$$ hence the result with
$1/C_{LP}=\alpha \, \left(1 - \frac{b \, \mu(U^c)}{\theta \, \alpha \, \mu(U)}\right)/(1+\lambda)$
since $\theta=b/\alpha$.

If $U$ does not contain the full level set $\{V \leq b/\alpha\}$ but only $C'$, the only
difference is that we cannot divide by $\theta$, hence the result.
\end{proof}

\begin{remark}
The conditions on $U$ are not really difficult to check in practice. We have included the first
situation because it covers cases where a bounded Lyapunov function exists, hence we cannot assume
in general that $U$ contains some level set. The second case is the usual one on euclidean spaces
when $V$ goes to infinity at infinity, so that we may always choose $U$ as a regular neighborhood
of a level set of $V$.

One may think that the constant $C_{LP}$ we have just obtained is a disaster. In particular,
contrary to the Meyn-Tweedie approach, the exponential rate given by $C_{LP}$ does not only depend
on $\alpha$ but also on $b,C,V$. But recall that in Meyn-Tweedie approach the non explicit
constant in front of the geometric rate depends on all these quantities (while we here have an
explicit $\int f^2 W d\mu$). In addition to the stronger type convergence ($\L^2$ type), one
advantage of Theorem \ref{thmlyapgap} is perhaps to furnish explicit (though disastrous)
constants.
\end{remark}
\medskip

%%%%%%%%%%%%%%%%%%%%%%%%%%%%%%%%%%%%%%%
%%%%%%%%%%%%%%%%%%%%%%%%%%%%%%%%%%%%%%%

\subsection{\bf A general sufficient condition for a Poincar\'e inequality.}\label{subspoincare}

As we previously said, there are some situations for which a tractable criterion for Poincar\'e's inequality is
known.

The most studied case is of course the euclidean space equipped with an absolutely continuous
measure $\mu(dx)=e^{-2F} dx$ and the usual $\Gamma(f)=|\nabla f|^2$. Dimension one is the only one
for which exists a general necessary and sufficient condition (Muckenhoupt criterion, see
\cite{logsob} Thm 6.2.2). A more tractable sufficient condition can be deduced (see \cite{logsob}
Thm 6.4.3) and can be extended to all dimensions using some isometric correspondence between
Fokker-Planck and Schr\"{o}dinger equations, namely $|\nabla F|^2(x) - \Delta F (x) \geq b > 0$
for all $|x|$ large enough (for a detailed discussion of the spectral theory of these operators
see \cite{HN} in particular Proposition 3.1). Actually this condition can be extended to $\mu =
e^{-2F} \nu$ if $\nu$ satisfies some log-Sobolev inequality, see \cite{GWu00} (as explained in
\cite{Cat5} Prop 4.4).

We shall see now that these conditions actually are of Lyapunov type, hence can be extended to a
very general setting.

\begin{lemma}\label{lemperturb}
Let $F$ be a nice enough function. Then if $V=e^{aF}$, $$LV \, - \, \Gamma(F,V) \, = \, a V \,
\left(LF \, + \, (\frac a2 \, - \, 1) \, \Gamma(F)\right) \, .$$
\end{lemma}
The proof is immediate using \eqref{eqcarre}, \eqref{eqderive} and \eqref{eqchainrule}. We may
thus deduce

\begin{theorem}\label{thmsuffgap}
Let $\nu$ be a ($\sigma$-finite positive measure) and $L$ be $\nu$ symmetric. Let $F \in D(L)$ be
non-negative and such that $\mu = (1/Z_F) \, e^{-2F} \nu$ is a probability measure for some
normalizing constant $Z_F$. For $0<a<2$ define
$$H_a=LF \, + \, (\frac a2 \, - \, 1) \, \Gamma(F)$$ and for $\alpha >0$, $C(a,\alpha)=\{H_a \geq
- (\alpha/a)\}$.

Assume that for some $a$ and some $\alpha$, $H_a$ is bounded above on $C(a,\alpha)$.

Assume in addition that for $\varepsilon > 0$ small enough one can find a large subset $U
\supseteq C(a,\alpha)$ with $\mu(U)\geq 1- \varepsilon$ such that $F$ is bounded on $U$, and $\mu$
satisfies some local Poincar\'e inequality on $U$.

Then $\mu$ satisfies some Poincar\'e inequality.
\end{theorem}
\begin{proof}
Recall that the operator  $L_Ff =Lf - \Gamma(F,f)$ is $\mu$ symmetric. According to Lemma
\ref{lemperturb}, if $V=e^{aF}$, $L_FV \leq - \alpha V$ outside $C(a,\alpha)$. But $H_a$ and $V$
being bounded on $C(a,\alpha)$, one can find some $b$ such that $V$ is a Lyapunov function. We may
thus apply Theorem \ref{thmlyapgap} which tells us that $\mu$ satisfies a Lyapunov-Poincar\'e
inequality. Since we are in the symmetric case, we may conclude thanks to Corollary
\ref{corpoinc}.
\end{proof}

We defer to Section 4 further results, applications and comments of this Theorem.
\medskip

%%%%%%%%%%%%%%%%%%%%%%%%%%%%%%%%%%%%%%%%%%%%%%%%%%%%%%%
%%%%%%%%%%%%%%%%%%%%%%%%%%%%%%%%%%%%%%%%%%%%%%%%%%%%%%%

\subsection{\bf Weak Lyapunov-Poincar\'e inequalities and weak Poincar\'e inequalities.}\label{subswpoincare}

We shall conclude this section by extending the two previous subsections to the more general weak
framework. We start with the following extension of Theorem \ref{thmlyapgap}

\begin{theorem}\label{thmwlyapgap}
Assume that there exists a $2 \phi$-Lyapunov function $V$ i.e. $LV \leq - 2 \phi(V) + b \BBone_C$
for some set $C$ (non necessarily petite). Recall that $\phi(u) \geq R >0$.

Assume that one can find a (large) set $U$ such that $\mu$ satisfies a local Poincar\'e inequality
on $U$.

Assume in addition that
\begin{enumerate}
\item[(1)] either $U$ contains $C'=C \cap \{\phi(V) \leq b\}$ and $R \, \mu(U) > b \mu(U^c)$ ,
\item[(2)] or $U$ contains $\{\phi(V) \leq b\}$ , $\phi$ is increasing and $\phi(b) \, \mu(U)> b
\, \mu(U^c)$.
\end{enumerate}
Then for $\lambda=(b \kappa_U - 1)_+$ and $W=V+\lambda$, $\mu$ satisfies a
(W)-weak-Lyapunov-Poincar\'e inequality, i.e. for all $f$ with $\int f d\mu =0$ and all $s>0$,
$$\int f^2 W d\mu \, \leq \, C_w \, \beta_W(s) \, \left(\int
\, \left(W \, \Gamma(f) \, - \, f^2 \, LW\right) \, d\mu\right) \, + \, s \, \parallel
f\parallel_\infty^2$$ with $\beta_W(s) = \inf \{u \, ; \, \int_{V > u \phi(V)} \, V d\mu \leq
s\}$, and where $C_w$ is given in the two corresponding cases by
\begin{enumerate}
 \item[(1)] $1/C_w= \left(1 - \frac{b \, \mu(U^c)}{R \,
\mu(U)}\right)/(1+\lambda)$, \item[(2)] or  $1/C_w= \left(1 - \frac{b \mu(U^c)}{\phi(b)
\mu(U)}\right)/(1+\lambda)$.
\end{enumerate}
\end{theorem}
\begin{proof}
Looking at the proof of Theorem \ref{thmlyapgap} we immediately see that, if $V$ is a $2
\phi$-Lyapunov function (recall definition \ref{defLyap}), then we may replace $C$ by $C'=C \cap
\{\phi(V) \leq b\}$. In the first situation we obtain as in the proof of Theorem \ref{thmlyapgap}
$$ \int_{C'}  \, f^2 \, d\mu \leq   \kappa_U \, \int \, \Gamma(f) \, d\mu \, + \,
(\mu(U^c)/R \, \mu(U)) \, \left(\int \, f^2 \, \phi(V) \, d\mu\right)\, ,$$ so that $$\int \,
\left(W \, \Gamma(f) \, - \, f^2 \, LW\right) \, d\mu \, \geq  \, \left(1 - \frac{b \, \mu(U^c)}{R
\, \mu(U)}\right) \, \int f^2 \, \phi(V) \, d\mu \, .$$ In the second case we may replace $R$ by
$\phi(b)$. It remains to note that $$ \int f^2 V d\mu \leq  u \, \int_{V \leq u \phi(V)} f^2
\phi(V) d\mu + \parallel f\parallel_\infty^2 \, \left(\int _{V > u \phi(V)} V d\mu\right) $$ for
all $u>0$.
\end{proof}

\begin{remark}
It is difficult to compare in full generality the previous weak Poincar\'e inequality with the one
obtained in Theorem \ref{thmMT}. More precisely, the previous result furnishes some decay for the
variance (as Theorem \ref{thmMT}) but the rate explicitly depends on $V$ (while $V$ only appears
through the constants in Theorem \ref{thmMT}. We shall thus make a more accurate comparison on
examples later on.

It is however worthwhile noticing that, in the first case, we do not need to impose any condition
on $\phi$ except that $\phi$ is bounded below by some positive constant.

Also remark that Theorem 3.1 in \cite{r-w} establishes a weak Poincar\'e inequality assuming that
one can find an exhausting sequence of sets $U_n$ such that $\mu$ satisfies a local Poincar\'e
inequality on each $U_n$. Here we only need one set $U$ (but large enough). Actually in the
examples we have in mind the assumption in \cite{r-w} is satisfied, but we shall see that we can
improve upon the function $\beta_W$.
\end{remark}

We shall now extend Theorem \ref{thmsuffgap} to the weak context.

\begin{corollary}\label{corsuffwgap}
Let $\nu$ be a ($\sigma$-finite positive measure) and $L$ be $\nu$ symmetric. Let $F \in D(L)$ be
non-negative and such that $\mu = (1/Z_F) \, e^{-2F} \nu$ is a probability measure for some
normalizing constant $Z_F$. We assume in addition that there exists $p<2$ such that $\int e^{-pF}
d\nu = c_p < +\infty$.

Let $\eta$ be a non-increasing function such that $u \, \eta(\log(u))$ is bounded from below by a
positive constant. For $0<a<2$ define $H_a=LF \, + \, (\frac a2 \, - \, 1) \, \Gamma(F)$ and
$C(a)=\{H_a \geq - \eta(F) \}$.

Assume that for some $0<a<2-p$ , $H_a$ is bounded above on $C(a)$. Assume in addition that for
$\varepsilon > 0$ small enough one can find a large subset $U \supseteq C(a)$ with $\mu(U)\geq 1-
\varepsilon$ and such that $F$ is bounded on $U$, and $\mu$ satisfies a local Poincar\'e
inequality on $U$.

Then $\mu$ satisfies  a weak Lyapunov-Poincar\'e, with $W=e^{aF}+\lambda$ (for some positive
$\lambda$),  inequality with \begin{equation}\label{bw}\beta_W(s)\, = \, \frac{2}{\left(a \,
\eta\left(\frac{\log(c_p/s )}{2-a-p}\right)\right)}\end{equation} hence for $\int f d\mu=0$, $\int
(P_t^* f)^2 \, d\mu \, \leq \, \int (P_t^* f)^2 \, W \, d\mu \, \leq \, \xi(t) \,
\parallel f\parallel^2_\infty$ with $$\xi(t)=2 \inf \{r>0 \, ; \, - C_w \, \beta_W(r) \, \log(r) \, \leq \,
t\} \, .$$ Finally $\mu$ satisfies a weak Poincar\'e inequality with $$\beta(s) = s \, \inf_{u>0}
\xi^{-1}(u \, e^{(1-u/s)}) \, \textrm{ with } \, \xi^{-1}(a) := \inf \{b>0 \, , \, \xi(b)\leq a\}
\, .$$
\end{corollary}

We easily remark that $\beta$ and $\beta_W$ are of the same order and change only through constants.

\begin{proof}
With our hypotheses, for $0<a<2$, $e^{aF}$ is a $2 \phi$-Lyapunov function for $\phi(u)=\frac 12
\, au \, \eta(\log(u)/a)$. Recall that we do not need here $\phi$ to be increasing nor concave. We
may thus apply Theorem \ref{thmwlyapgap} yielding some weak Lyapunov-Poincar\'e inequality for
$\mu$.

We shall describe the function $\beta_W$. Recall that
\begin{eqnarray*}
\beta_W(s) & = & \inf \{u \, ; \, \int_{V > u \phi(V)} \, V d\mu \leq s\}\\ & = & \inf \{u \, ; \,
\int_{2 > au \eta(F)} \, e^{(a-2)F} d\nu \leq s\}\\ & = & \inf \{u \, ; \, \int_{F >
\eta^{-1}(2/au)} \, e^{(a-2)F} d\nu \leq s\} \, .
\end{eqnarray*}
But if $2-a=p+m$,
\begin{eqnarray*}
\int_{F > \eta^{-1}(2/au)} \, e^{(a-2)F} d\nu & \leq & \int \, e^{-pF} \, e^{- m \eta^{-1}(2/au)}
\, d\nu
\end{eqnarray*}
from which we deduce that $\beta_W(s) \leq \, 2/(a \eta(\frac 1m \, \log(c_p/s)))$.
\smallskip

Using Lemma \ref{lemmecle} we deduce as usual (see e.g. the proof of Theorem 2.1 in \cite{r-w})
that $\int (P_t^* f)^2 \, W \, d\mu \, \leq \, \xi(t) \, \parallel f\parallel^2_\infty$ with
\begin{equation}\label{eqvitesse}
\xi(t)=2 \inf \{r>0 \, ; \, - C_w \, \beta_W(r) \, \log(r) \, \leq \, 2t\}
\end{equation}
for $\int f d\mu=0$, $W$ and $C_w$ being as in the previous Theorem.
\end{proof}
\smallskip

\begin{remark}
In view of Theorem \ref{thmMT} it is interesting to replace the $\L^\infty$ norm above by $\L^p$
norms, with $p>2$. In the case of usual weak Poincar\'e inequalities it is known that we may
replace the $\L^\infty$ norm by a $\L^p$ norm just changing the $\beta$ into $\beta_p(s)=c \,
\beta (c' s^q)$ for some constants $c$ and $c'$, and $1/p+1/q=1$ (see e.g. \cite{Zitt} Theorem 29
for a more general result). But the proof in \cite{Zitt} (inspired by \cite{CGG} Theorem 3.8) lies
on a Capacity-Measure characterization of these inequalities introduced in \cite{BCR2}.

The situation here is more complicated and a direct modification of the weak-Lyapunov-Poincar\'e
inequality seems to be difficult. However, since we are interested in the rate of convergence to
the equilibrium, we may mimic the truncation argument in \cite{CGG}. Namely, let $f$ be such that
$\int f d\mu =0$, denote by $f_K=f\wedge K \vee -K$ and $m_K=\int f_K d\mu$, then if a
weak-Lyapunov-Poincar\'e inequality holds we get for all $p>1$,
\begin{eqnarray*}
\int (P_t^*f)^2 d\mu  \, \, \leq & 2 \, \left(\int (P_t^* (f_K-m_K))^2 W d\mu + \int (P_t^*
(f-f_K+m_K))^2 d\mu \right) \, \\ \leq & 2 \, \left(\int (P_t^* (f_K-m_K))^2 W d\mu + \int
(f-f_K+m_K)^2 d\mu \right) \\ \leq & 2 \, \xi(t) \, K^2 + 4 \, \int_{|f|>K} \, (|f|-K)^2 \, d\mu
\, + \, 4 m_K^2  \\ \leq & 2 \, \xi(t) \, K^2 + 4 \, \int_{|f|>K} \, (|f|-K)^2 \, W d\mu \, + \, 4
\, \left(\int_{|f|>K} ||f| - K| d\mu\right)^2 \\ \leq & 2 \, \xi(t) \, K^2 + 8 \, \int_{|f|>K} \,
(|f|-K)^2 \, d\mu \\ \leq & 2 \xi(t) \, K^2 + 8 \left(\int |f|^{2p} d\mu\right) \, K^{-2p/q} \, .
\end{eqnarray*}
Now optimizing in $K$ furnishes
\begin{equation}\label{eqdectrunc}
\int (P_t^*f)^2 d\mu \, \leq \, C \, \xi^{1/q}(t) \, \left(\int |f|^{2p} d\mu\right)^{1/p}
\end{equation}
which is quite the result in Theorem \ref{thmMT}, but with explicit constants.
\end{remark}

\begin{remark}

It is perhaps more natural to try to obtain directly a weak Poincar\'e inequality starting from
the existence of a $\phi$-Lyapunov function as follows.

For $\int f d\mu=0$, we have $$\int f^2 d\mu \leq \int   \, \frac{-LV}{\phi(V)} f^2 \, d\mu + \int
f^2 \frac{b}{\phi(V)} \, \BBone_C \, d\mu \, .$$ We know how to manage the second term if a local
Poincar\'e inequality holds, hence we focus on the first term in the right hand side of the
previous inequality.

Assume that $L$ is $\mu$-symmetric. Integrating by parts we get
\begin{eqnarray*}
\int  \, \frac{-LV}{\phi(V)}  f^2 \, d\mu & = & \int \left(\frac{f \, \Gamma(f,V)}{\phi(V)} \, -
\, \frac{f^2 \phi'(V) \Gamma(V)}{2 \phi^2(V)}\right) d\mu
\end{eqnarray*}
but thanks to our hypotheses $$\frac{f \, \Gamma(f,V)}{\phi(V)} \leq \frac a2 \, \Gamma(f) +
\frac{1}{2a} \, \frac{f^2 \, \Gamma(V)}{\phi^2(V)}$$ for all $a>0$ so that
\begin{eqnarray*}
\int  \, \frac{-LV}{\phi(V)}  f^2 \, d\mu & \leq & \int \frac a2 \, \Gamma(f) d\mu \, + \, \int \,
\left(\frac{f^2 \, \Gamma(V)}{\phi^2(V)}\right)\left(\frac{1}{2a} - \phi'(V)\right) d\mu
\end{eqnarray*}
Unless $\phi$ is linear, $\liminf \phi'=0$ at infinity, so that we get an extra term that cannot
be controlled. Of course if $\phi$ is linear, we may choose $a$ for this term to vanish
identically, and so obtain another proof of the Poincar\'e inequality for $\mu$ (with more easily
calculable constants).
\end{remark}
\bigskip

%%%%%%%%%%%%%%%%%%%%%%%%%%%%%%%%%%%%%%%%%%
%%%%%% exemples...
%%%%%%
%%%%%%%%%%%%%%%%%%%%%%%%%%%%%%%%%%%%%%%%

\section{\bf Examples.}\label{secex1}

Due to the local Poincar\'e property, the most natural framework is the euclidean space $\R^d$. It
will be our underlying space in all examples, but in many cases results extend to a Riemanian
manifold as well.

\subsection{\bf General weighted Poincar\'e inequalities.}

Let $F$ be a smooth enough non-negative function such that $\mu = (1/Z_F) \, e^{-2F} dx$ is a
probability measure. We may also assume that $F(x) \to +\infty$ as $x \to \infty$, so that the
level sets of $F$ are compact. If
\begin{itemize}
\item either $\Delta \, F - \, |\nabla F|^2$ is bounded from above, \item or $\int |\nabla F|^2
d\mu \, < +\infty$ ,\end{itemize} it is known than one can build a conservative (i.e. non
exploding) $\mu$ symmetric diffusion process with generator $L_F = \frac 12 \, \Delta - \nabla
F.\nabla$. We shall assume for simplicity that the first condition holds.

Assume in addition that $$\liminf_{|x| \to +\infty} \left(|\nabla F|^2 \, - \, \Delta \, F \right)
= \alpha > 0 \, .$$ We may thus apply Theorem \ref{thmsuffgap}, with $L=\frac 12 \, \Delta$, $\nu$
the Lebesgue measure (which is known to satisfy a Poincar\'e inequality on euclidean balls of
radius $R$ with $C_P=C R^2$ , $C$ being universal, and for $\Gamma (f)=|\nabla f|^2$), $U$ a large
enough ball, $a=1$. Indeed, since $\nu$ satisfies a (true) Poincar\'e inequality on $U$, $\mu$
which is a log-bounded perturbation of $\nu$ on $U$ also satisfies a Poincar\'e inequality on $U$,
hence a local one (since the energy on $U$ is smaller than the one on the full $E$). This yields
\begin{corollary}\label{corweight1}
If $F$ is a $C^2$ non-negative function such that, $F(x) \to +\infty$ as $x \to \infty$, $\int
e^{-2F} dx < + \infty$ and \begin{itemize} \item $\Delta \, F - \, |\nabla F|^2$ is bounded from
above, and , \item $\liminf_{|x| \to +\infty} \left(|\nabla F|^2 \, - \, \Delta \, F \right) =
\alpha > 0 \, .$ \end{itemize} Then the following (weighted) Poincar\'e inequality holds for all
$f$ smooth enough and some $C_P$, $$\int f^2 \, e^{-2F} \, dx \, \leq \, C_P \, \int |\nabla f|^2
\, e^{-2F} \, dx \, + \, \frac{\left(\int f \, e^{-2F} \, dx\right)^2}{\left(\int e^{-2F}
dx\right)} \, .$$
\end{corollary}

This corollary immediately extends to uniformly elliptic operators in divergence form. The
degenerate case is more intricate. Indeed, according to results by Jerison, Franchi, Lu
(\cite{Jer,Fra,Lu}) a Poincar\'e inequality holds on small metric balls for more general operators
of locally subelliptic type. Let us describe the framework we are interested in.

Let $X_1,...,X_m$ be $C^\infty$ vector fields defined on $\R^d$. We shall assume for simplicity
that they are bounded with all bounded derivatives. We shall make the following H\"{o}rmander type
assumption : \begin{assumption}\label{asshormander} there exists $N \in \N^*$ and $c>0$ such that
for all $x$ and all $\xi \in R^d$,
 $\sum_Y \, \langle Y(x),\xi\rangle^2 \, \geq \, c |\xi|^2$, where the sum is taken
over all Lie brackets $Y=[X_{i_1},[...X_{i_k}]]$ of length less than or equal to $N$.
\end{assumption}

This assumption is enough for ensuring that the natural associated subriemanian metric $\rho$ is
locally equivalent to the usual one (see e.g. \cite{Fra} Theorem 2.3). According e.g. to Theorem C
in \cite{Lu} (a similar result was first obtained by Jerison), the Lebesgue measure $dx=\nu$
satisfies a Poincar\'e inequality on small metric balls $B_\rho(y,s)$ for $s$ small enough and
$\Gamma(f)=\sum_{i=1}^m \, |X_i f|^2$. But here we want some local Poincar\'e inequality on some
large set.

If we replace the euclidean space by a connected unimodular Lie group with polynomial volume
growth equipped with left invariant vector fields $X_1,...,X_m$ generating the Lie algebra of $E$,
then it is known that a Poincar\'e inequality holds for all metric balls (the result is due to
Varopoulos and we refer to \cite{SC94} p.275 for explanations). But in the euclidean case we can
show that Lebesgue measure satisfies some local Poincar\'e inequality on euclidean balls centered
at the origin.

Indeed let $|.|$ stands for the euclidean norm. Recall that there exist $R$ and $r$ such that
$$\{|x|\leq r\} \, \subset B_\rho(0,s) \, \subset \{|x|\leq R\} \, .$$ If $\int_{|x|\leq N} f dx =0$,
then for all $a$ it holds
\begin{eqnarray*}
\int_{|x|\leq N} f^2(x) \,  dx & = & \int_{|x|\leq r} f^2(Nx/r) \, (\frac Nr)^d \, dx \\ & \leq &
\int_{|x|\leq r} (f(Nx/r)-a)^2 \, (\frac Nr)^d \, dx \\ & \leq & \int_{B_\rho(0,s)} (f(Nx/r)-a)^2
\, (\frac Nr)^d \, dx \, ,
\end{eqnarray*}
so that if we choose $a= (\int_{B_\rho(0,s)} f(Nx/r) \, dx)/|B_\rho(0,s)|$ (where $|U|$ denotes
the Lebesgue volume of $U$) we may use the Poincar\'e inequality in the metric ball, and obtain
(denoting by $g(x)=f(Nx/r)$)
\begin{eqnarray*}
\int_{|x|\leq N} f^2(x) \,  dx & \leq & C \, (\frac Nr)^d \, \int_{B_\rho(0,s)} \sum_{i=1}^m |X_i
g|^2(x) \,dx \\ & \leq & C \, \left(\frac Nr\right)^{d+2} \, \int_{B_\rho(0,s)} \sum_{i=1}^m |X_i
f|^2(Nx/r) \,dx \\ & \leq & C \, \left(\frac Nr\right)^{d+2} \, \int_{|x|\leq R} \sum_{i=1}^m |X_i
f|^2(Nx/r) \,dx \\ & \leq & C \, \left(\frac Nr\right)^2 \, \int_{|x|\leq (RN/r)} \sum_{i=1}^m
|X_i f|^2(x) \,dx \, ,
\end{eqnarray*}
Now, since $F$ is locally bounded, it is straightforward to show that $\mu$ satisfies a local
Poincar\'e inequality on $\{|x|\leq N\}$ with $\kappa_N = C \left(\frac Nr\right)^2 \, e^{4
\sup_{|x|\leq (RN/r)}F(x)}$.
\smallskip

If we define a new vector field as $X_0 = \frac 12 \, \sum_{i=1}^m \, div X_i .X_i$ , then $dx$ is
symmetric for the generator $L = \frac 12 \, \sum_{i=1}^m \, X_i^2 + X_0$ and $\mu$ is symmetric
for the generator $L_F = \frac 12 \, \sum_{i=1}^m \, X_i^2 + X_0 - \sum_{i=1}^m \, X_i F . X_i$
written in H\"{o}rmander form. Hence the following generalizes Corollary \ref{corweight1}

\begin{corollary}\label{corweight2}
Assume that Assumption \ref{asshormander} is fulfilled, and let $L = \frac 12 \, \sum_{i=1}^m \,
X_i^2 + X_0$ be as above. If $F$ is a $C^2$ non-negative function such that, $F(x) \to +\infty$ as
$x \to \infty$, $\int e^{-2F} dx < + \infty$ and
\begin{itemize} \item $LF - 1/2 \sum_{i=1}^m \, |X_i F|^2 $ is bounded from above, and , \item
$\liminf_{|x| \to +\infty} \left( 1/2 \sum_{i=1}^m \, |X_i F|^2 - LF\right) = \alpha > 0 \, .$
\end{itemize} Then the following (weighted) Poincar\'e inequality holds for all $f$ smooth enough
and some $C_P$, $$\int f^2 \, e^{-2F} \, dx \, \leq \, C_P \, \int \, \sum_{i=1}^m \, |X_i f|^2 \,
e^{-2F} \, dx \, + \, \frac{\left(\int f \, e^{-2F} \, dx\right)^2}{\left(\int e^{-2F} dx\right)}
\, .$$
\end{corollary}

\begin{remark}
The choice $V=e^{aF}$ is not necessarily the best possible. Indeed one wants to get the smallest
possible Lyapunov function. For example if $F(x)=|x|^2$ (i.e. the gaussian case) one can choose
$V(x)=1+a|x|^2$ for $a>0$. This is related to some sufficient condition for the Gross logarithmic
Sobolev inequality to hold (see \cite{Cat5}). In the same way, if $F(x)=|x|^p$ (at least away from
0 for $F$ to be smooth) for some $2>p\geq 1$, it is easy to see that $V(x)=\exp \left(a
|x|^{2-p}\right)$ is a Lyapunov function (at least for a good choice of $a$), and of course $2-p <
p$ when $p>1$, so that this choice is better than $e^{aF}$. These laws of exponent $1<p<2$ are the
generic examples of laws satisfying interpolating inequalities (called $F$-Sobolev inequalities
see \cite{BCR1}, take care that this $F$ is not the potential). It clearly suggests that the best
possible choice for the Lyapunov function is connected with the $F$-Sobolev inequality satisfied
by $\mu$.
\end{remark}

%%%%%%%%%%%%%%%%%%%%%%%%%%%%%%%
%%%%%%%%%%%%%%%%%%%%%%%%%%%%%%%

\subsection{\bf General weighted weak Poincar\'e inequalities.}

In this subsection we shall compare various weak Poincar\'e inequalities obtained in \cite{r-w},
\cite{BCR2}, Theorem \ref{thmMT} and Corollary \ref{corsuffwgap}; as well as the various rates of
convergence to equilibrium. The framework is the same as in the previous subsection.
\smallskip

\subsubsection{\bf Sub-exponential laws.}

We consider here for $0<p<1$ the measures $\mu_p(dx)=C_p \, e^{-2|x|^p} dx$ where $C_p$ is a
normalizing constant and $|.|$ denotes the euclidean norm. It is shown in \cite{BCR2} that if
$d=1$, $\mu_p$ satisfies a weak Poincar\'e inequality with $\beta_p(s)=d_p \, \log^{(2/p)-2}(2/s)$
this function being sharp. Note that the previous result does not extend to higher dimensions via
the tensorization result (Theorem 3.1 in \cite{BCR2}). In any dimension however, \cite{r-w}
furnishes $\beta_p(s)=d_p \, \log^{(4(1-p)/p}(2/s)$. Note that for $d=1$ the result in \cite{BCR2}
improves on the one in \cite{r-w}.

These bounds furnish a sub-exponential decay $$\int (P_t^* f)^2 \, d\mu \, \leq \, c_1 \, e^{-c_2
\, t^\delta} \,
\parallel f\parallel^2_\infty$$ for any $\mu$ stationary semi-group, with $\delta=p/(2-p)$ if
$d=1$ and $\delta=p/(4-3p)$ for any $d$.
\smallskip

But sub-exponential laws enter the framework of subsection \ref{subswpoincare} with $V=e^{a|x|^p}$
, $\eta(u)= c \, u^{2(1-\frac 1p)}$ hence $\beta_W(s)= C \, \log^{(2/p)-2}(c/s)$ for some
constants $c$ and $C$. Note that we recover the right exponent $(2/p)-2$ for $\beta_W$, hence the
right sub-exponential decay in any dimension. Up to the constants, we also recover, thanks to
Theorem 2.3 in \cite{r-w}, that $\beta$ behaves like $\beta_W$.

Also note that in this case the rate given by Theorem \ref{thmMT} is again $\psi(t) = c_1 \, e^{-
c_2 \, t^{p/(2-p)}}$.
\smallskip

These results extend to any $F$ going to infinity at infinity and satisfying $$(1-a/2)|\nabla F|^2
- \Delta F \geq c \, F^{2(1-\frac 1p)}$$ at infinity, generalizing to the weak Poincar\'e
framework similar results for super-Poincar\'e inequalities (see \cite{BCR1} and \cite{BCR3}).

\subsubsection{\bf Heavy tails laws.}

Let us deal now with measures $\mu_p(dx)=C_p(1+|x|)^{-(d+p)}$ where $p>0$, $C_p$ is a normalizing
constant,  and $|.|$ denotes once again the usual euclidian norm. The sharp result in dimension 1
has been given in \cite{BCR2} with $\beta_p(s)=d_ps^{-2/p}$, but cannot be extended to higher
dimensions. R\"ockner-Wang \cite{r-w} furnishes in any dimension $\beta_p(s)=cs^{-\tau}$ where
$\tau=\min\{(d+p+2)/p,(4p+4+2d)/(p^2-4-2d-2p)_+\}$. This result is not sharp in dimension one but
enables to quantify the polynomial decay of the variance in any dimension.

Once again, we may use the results of section  \ref{subswpoincare} with $V(x)=(1+|x|)^{a(d+p)/2}$,
so that  $F(x)={d+p\over2}\log(1+|x|)$ and $\eta(u)=C(p,d)e^{-4u/(p+d)}$. Use now (\ref{bw}) to
get that $\beta_W(s)=C(p',d)s^{2\over p'}$ for any $p'<p$ (and $C(p',d)\to\infty$ as $p'\to p$).
This result enables us to be nearly optimal in any dimension and thus improves on the result of
\cite{r-w}. Note that once again, results of \cite{Forro,DFG} would give, via Theorem \ref{thmMT}
the same result, but without explicit constants.
\medskip

\subsection{\bf Drift conditions for diffusion processes.}

Consider a $d$ dimensional diffusion process
\begin{equation}\label{eqdiffusion}
dX_t \, = \, \sigma(X_t) dB_t + \beta(X_t) dt \, .
\end{equation}
We assume that the (matrix) $\sigma$ has smooth and bounded entries, and is either uniformly
elliptic or hypoelliptic in the sense of Assumption \ref{asshormander}. We also assume the
following drift condition
\begin{equation}\label{eqdriftcond}
\textrm{ there exist $M$ and $r>0$ such that for all $|x|\geq M$ , } \langle \beta(x) , x \rangle
\, \leq \, - r |x| \, .
\end{equation}
We also assume that the diffusion has an unique invariant probability measure $d\mu=e^{F} dx$.
This is automatically satisfied if $\sigma$ is uniformly elliptic and \eqref{eqdriftcond} holds
(see \cite{DFG} Proposition 4.1).
\medskip

Consider a smooth function $V$ which coincides with $e^{a |x|}$ outside the ball of radius $M$,
$|x|$ denoting the euclidean distance. Then on this set $$LV(x)=a \langle \beta(x),
\frac{x}{|x|}\rangle + a^2 \eta(x)$$ where $\eta(x) \to 0$ as $|x| \to +\infty$. Hence, according
to \eqref{eqdriftcond} for all $a$, $V$ is a Lyapunov function (but $C$ and $b$ depend on $a$).

We may thus apply Theorem \ref{thmlyapgap} (thanks to the local Poincar\'e property discussed in
the previous subsection) and get that for any density of probability $h$, $$\int |P^*_t h - 1|^2
\,  d\mu \, \leq \, e^{-\delta_a t} \, \int (h-1)^2 \, e^{a |x|} \, d\mu \, .$$ Indeed we know
that $\mu$ satisfies a $(V+\lambda)$-Lyapunov-Poincar\'e inequality, hence apply Proposition
\ref{proptrou} and then replace $W$ by $1+\lambda$ in the left hand side, and $(V+\lambda)$ by
$(1+\lambda)V$ in the right hand side. Hence we get an exponential convergence for initial
densities in $\L^2(e^{a|x|} \mu)$ for some $a>0$.

Remark that if $\sigma=Id$ and $\beta=-\nabla F$, $d\mu=e^{-2F} dx$ and \eqref{eqdriftcond} which
reads $$\textrm{ there exist $M$ and $r>0$ such that for all $|x|\geq M$ , } \langle \nabla F(x) ,
x \rangle \, \geq \,  r |x| \, ,$$ thus implies the Poincar\'e inequality.

We may now complete the picture in the sub-exponential case (the polynomial case being handled
similarly),  namely we assume

\begin{equation}\label{eqdriftcondsub}
\textrm{ there exist $0<p<1$, $M$ and $r>0$ such that for all $|x|\geq M$,  , } \langle \beta(x) , x \rangle
\, \leq \, - r |x|^{1-p} \, .
\end{equation}
One may then show as in \cite{DFG} that for sufficiently small $a$, $V(x)=e^{a|x|^{1-p}}$ is a
$\phi$-Lyapunov  function with $\phi(v)=v \, \log(v)^{-2{p\over 1-p}}$ and get via the use of
theorem \ref{thmwlyapgap} as in the preceding paragraph a weak Lyapunov-Poincar\'e inequality with
$W(x)=V(x)+\lambda$ and $\beta_W(s)=d_p\log(2/s)^{2p/(1+p)}$. It then implies that for any density
of probability $h$, $$\int |P^*_t h - 1|^2 \,  d\mu \, \leq \, C_{a,p}e^{-\delta_a t^{1-p\over
1+p}} \, \int (h-1)^2 \, e^{a |x|^{1-p}} \, d\mu \, .$$ Let us remark that for this diffusion
case, the use of Lyapunov function was already present in R\"ockner-Wang \cite[Th. 3.2 and 3.3]{r-w} to obtain weak Poincar\'e inequality. They however always propose as Lyapunov function the distance to
the origin, combined with local approximations, which is not optimal as seen in the previous
subsections. Remark however that as in \cite{DFG}, R\"ockner-Wang \cite{r-w} also considers the
case of Markov processes with jumps. We leave this for further research.
\bigskip

\section{\bf Entropy and weighted entropy.}\label{secentrop}

In all the previous sections we studied the behaviour of the Variance or some weighted Variance.
The only exception is Theorem \ref{thmMTlog} where we obtained the rate of convergence for
relative entropy. In many significant cases, for physical relevance, $\L^2$ bounds are too
demanding, so that it is of some interest to look at less demanding bounds.
\smallskip

Using Lemma \ref{lemmecle} the following Proposition is obtained exactly as Proposition
\ref{proptrou}, after stating the analogue of Definition \ref{defLP}
\begin{definition}\label{defLpsi}
Let $\Psi$ be a non-negative function such that $\Psi(1)=0$. We shall say that $\mu$ satisfies a
(W)-Lyapunov-$\Psi$-Sobolev inequality, if there exists $W\in D(L)$ with $W\geq 1$ and a constant
$C_{\Psi}$ such that for all nice non-negative $h$ with $\int h d\mu=1$,
$$\int \, \Psi(h) \, W \, d\mu \, \leq \, C_{\Psi} \, \int \, \left( \frac12 \, W \, \Psi''(h) \,
\Gamma(h) \, - \, \Psi(h) \, LW\right) \, d\mu \, .$$
\end{definition}

\begin{proposition}\label{proppsi}
Let $\Psi$ be a non-negative function such that $\Psi(1)=0$. The following statements are
equivalent
\begin{itemize} \item $\mu$ satisfies a (W)-Lyapunov-$\Psi$-Sobolev inequality, \item $\int
\Psi(P_t^*h) \, W \, d\mu \, \leq \, e^{- (t/C_{\Psi})} \, \int \Psi(h) W d\mu$ for all
non-negative $h$ with $\int h d\mu = 1$.
\end{itemize}
\end{proposition}

Since the goal of this section is to deal with densities of probability $h$ with very few moments
(in particular not in $\L^2$), we shall not discuss the analogous weak versions of these
inequalities. The interested reader will easily derive the corresponding results.
\smallskip

Note that for Definition \ref{defLpsi} to be interesting, we do certainly have to assume that
$\Psi''(u)>0$ for all $u$. This is a big difference with the (homogeneous) $F$-Sobolev
inequalities studied in \cite{BCR1} where $F$ often vanishes on some neighborhood of $0$.

Indeed if we want to mimic what we have done in Theorem \ref{thmlyapgap}, we have to introduce
some local version of some new $\Psi$-Sobolev inequality, replacing the local Poincar\'e
inequality. Instead of looking at such a complete theory, we shall focus on a typical example
which will give the flavor of the results one can obtain. The first remark is, see for instance
\cite{G1},  that the Lebesgue measure satisfies a logarithmic Sobolev inequalities on the interval
$I=[-R,R]$ with constant $8R^2/\pi^2$ which by tensorization holds also on the tensor product
$I^d$ with the same constant so that we obtain the equivalent of the local Poincar\'e inequality.

% Our first aim is to show that the
%Lebesgue measure satisfies ad-hoc inequalities on hypercubes.
%\medskip
%
%First we recall some definitions. Given $A\subset U$ we define $$Cap(A):= \inf \, \{Cap(A,\Omega);
%\, A \subset \Omega \, , \, |\Omega| \leq 1/2 \}$$ where $|\Omega|$ is the normalized Lebesgue
%measure of $\Omega$ and $$Cap(A,\Omega) := \inf \left\{ \int |\nabla f|^2 d\nu \, ; \, \BBone_A
%\leq f \leq \BBone_\Omega\right\}$$ where the infimum is taken over locally Lipschitz functions
%and where $\nu$ is the normalized Lebesgue measure.
%\smallskip
%
%Consider the interval $I=[-R,R]$ equipped with the normalized Lebesgue measure $d\nu=dx/2R$. If
%$A\subset I$ is non void, with the previous notations we may choose some $a\in A$ and $\omega \in
%\Omega^c$, so that $1 = f(a) - f(\omega) = \int_\omega^a \, f'(t) \, dt$. Hence, $$\int_{-R}^R \,
%|f'(t)|^2 dt \, \geq \, \int_{\omega}^a \, |f'(t)|^2 dt \, \geq \, \frac{1}{|a-\omega|} \, \geq \,
%\frac{1}{2R} \, ,$$ so that $$Cap(A) \, \geq \, \frac{1}{4R^2} $$ for all $A$.
%
%In particular $$\nu(A) \, \log\left(\frac{2}{\nu(A)}\right) \, \leq \, \log 2 \, R^2 \, Cap(A) \,
%.$$ In addition we know that $\nu$ satisfies on $I$ a Poincar\'e inequality with a constant $C_P
%\leq 2R^2$. Thus according to Theorem 26 in \cite{BCR1}, $\nu$ satisfies a log-Sobolev inequality
%with a constant $C_{LS} \leq C R^2$ where $C$ is universal ($C=60$ works). The log-Sobolev
%inequality being dimension free, it holds on the tensor product $I^d$ with the same constant. This
%result is of course known.
\smallskip

Now if $d\mu = e^{-2F} dx$ is a Probability measure, the normalized measure
$\bar{\mu}=\mu/\mu(I^d)$ also satisfies a log-Sobolev inequality on $I^d$ as soon as $F$ is
locally bounded.

But $u \mapsto u \log u$ is not everywhere non-negative so that we have to modify it.

First, since $\bar{\mu}$ also satisfies a Poincar\'e inequality on $I^d$, we may apply Lemma 17 in
\cite{BCR3} and obtain the following $G$-Sobolev inequality with $G(u)=(\log u - \log 4)_+$ and
some universal $C$ (all universal constants will be denoted by $C$ in the sequel)
\begin{equation}\label{eqFsobent}
\int_{I^d} \, f^2 \, G\left(\frac{f^2}{\int_{I^d} f^2 d\bar{\mu}}\right) \, d\bar{\mu} \, \leq \,
C \, (1+R^2) \,  \int_{I^d} \, |\nabla f|^2 \, d\bar{\mu} \, .
\end{equation}

Now consider $\Psi$ defined on $\R^+$ by
\begin{equation}\label{eqPsi}
\Psi(u)= (u-1)^2 \, \BBone_{u\leq 2} + \left(1 + (1-4\log 2)(u-2) + 4(u \log u - u -2 \log 2
+2)\right) \, \BBone_{u>2} \, ,
\end{equation}
so that $$\Psi''(u)=2 \BBone_{u\leq 2} + \frac{4}{u} \, \BBone_{u>2} \, ,$$ is everywhere
positive. $\Psi$ is non-negative and $\Psi(u)=0$ if and only if $u=1$. It is easy to see that $u
\mapsto \Psi(u)/u$ is non-decreasing on $[1,+\infty[$ and of course $\Psi$ behaves like $4 G$ at
infinity. Thus combining \eqref{eqFsobent} and Lemma 21 in \cite{BCR1} we obtain that for any nice
$g$ with $\int_{I^d} g^2 d\bar{\mu}=1$ it holds
\begin{equation}\label{eqPsisob}
\int_{I^d} \, \Psi(g^2) \, \BBone_{g^2 > 1} \, d\bar{\mu} \, \leq \, C \, (1+R^2) \, \int_{I^d} \,
|\nabla g|^2 \, d\bar{\mu} \, .
\end{equation}
We may thus state

\begin{theorem}\label{thmentropydecaywithlyap}
Let $\mu = e^{-2F} dx$ be a probability measure on $\R^d$ (supposed to be $L$ invariant)
satisfying a Poincar\'e inequality (on the whole $\R^d$) with constant $C_P$. Assume that there
exists a Lyapunov function $V$ i.e. $LV \leq - 2 \alpha V + b \BBone_C$ for some set $C$ (non
necessarily petite), such that either $C$ or the level sets of $V$ are compact.

Then $\mu$ satisfies a (W)-Lyapunov-$\Psi$-Sobolev inequality for $W=V+\lambda$ where $\lambda$ is
a large enough constant and $\Psi$ is defined in \eqref{eqPsi}.
\end{theorem}

\begin{remark}
According to Corollary \ref{corpoinc} and Theorem \ref{thmlyapgap}, if $L$ is $\mu$ symmetric, the
Poincar\'e inequality automatically holds here.
\end{remark}

\begin{proof}
Since we assumed that $C$ or the level sets of $V$ are compact, as for the proof of Theorem
\ref{thmlyapgap} what we have to do is to control $\int_{I^d} \Psi(h) d\mu$ for a large enough
$I^d$ and a non-negative $h$ such that $\int_{\R^d} \, h d\mu =1$. In the sequel we write $h=f^2$
(we may first assume that $f\geq \varepsilon > 0$ and then go to the limit if necessary).

First, applying Poincar\'e inequality we get
\begin{eqnarray*}
\int_{I^d} \, \Psi(h) \, \BBone_{h\leq 2} \, d\mu & = & \int_{I^d} \, (h-1)^2 \, \BBone_{h\leq 2}
\, d\mu \, = \,  \int_{I^d} \, (h\wedge 2 - 1)^2 \, \BBone_{h\leq 2} \, d\mu \\ & \leq &
\int_{\R^d} \, (h\wedge 2-1)^2 \, d\mu
\\ & \leq & C_P \, \int_{\R^d} |\nabla h|^2 \, \BBone_{h\leq 2} \, d\mu + \left(\int_{\R^d}
(h\wedge 2 - 1) \, d\mu\right)^2 \\ & \leq & C_P \, \int_{\R^d} |\nabla h|^2 \, \BBone_{h\leq 2}
\, d\mu + \left(\int_{\R^d} \left((h - 1) \BBone_{h < 2}+\BBone_{h\geq 2}\right) \, d\mu\right)^2\\
& \leq & C_P \, \int_{\R^d} |\nabla h|^2 \, \BBone_{h\leq 2} \, d\mu + \left(\int_{\R^d} \left((1
- h) \BBone_{h\geq 2}+\BBone_{h\geq 2}\right) \, d\mu\right)^2\\& \leq & C_P \, \int_{\R^d}
|\nabla h|^2 \, \BBone_{h\leq 2} \, d\mu + \left(\int_{\R^d} (2 - h) \, \BBone_{h\geq 2} \,
d\mu\right)^2 \\ & \leq & C_P \, \int_{\R^d} |\nabla h|^2 \, \BBone_{h\leq 2} \, d\mu +
\left(\int_{\R^d} (h - 2) \, \BBone_{h\geq 2} \, d\mu\right)^2 \\ & \leq & C_P \, \int_{\R^d}
|\nabla h|^2 \, \BBone_{h\leq 2} \, d\mu + \left(\int_{\R^d} h  \, \BBone_{h\geq 2} \,
d\mu\right)^2 \\ & \leq & C_P \, \int_{\R^d} |\nabla h|^2 \, \BBone_{h\leq 2} \, d\mu +
\int_{\R^d} h \, \BBone_{h\geq 2} \, d\mu \, ,
\end{eqnarray*}
since $\int_{\R^d} h \, \BBone_{h\geq 2} \, d\mu \leq 1$. Since $\mu$ satisfies a Poincar\'e
inequality, Remark 22 in \cite{BCR1} shows that $$\int f^2 \, \BBone_{f^2 \geq 2 \int f^2 d\mu} \,
d\mu \, \leq \, C \, \int |\nabla f|^2 d\mu \, ,$$ so that (recall that $\Psi''(u)=2 \BBone_{u\leq
2} + \frac{4}{u} \, \BBone_{u>2}$) we finally obtain for some constant $C$, $$\int_{I^d} \,
\Psi(h) \, \BBone_{h\leq 2} \, d\mu \, \leq \, C \, \int_{\R^d} \, \Psi''(h) \, |\nabla h|^2 \,
d\mu \, .$$ For the other part we have to be accurate with normalization in order to use
\eqref{eqPsisob}. Indeed the latter applies for normalized functions for the normalized measure on
$I^d$.

Let $m=\int_{I^d} (h\vee a) \, d\bar{\mu}$ for some $2>a>0$.

If $m\leq 1$ then $$\Psi(h) \, \BBone_{h>2} = \Psi(h\vee a) \, \BBone_{h>2}  \leq \Psi(h\vee a/m)
\, \BBone_{h>2} \leq \Psi(h\vee a/m) \, \BBone_{(h\vee a/m)>2}$$ so that we may apply
\eqref{eqPsisob} with $g=(h\vee a/m)^{\frac 12}$ (we can of course replace $\BBone_{g>1}$ by
$\BBone_{g>2}$). Of course $|\nabla g|^2$ is up to some constant (the normalization by $m$
disappears) equal to $\BBone_{h>a} (|\nabla h|^2/h)$ hence up to the constants to $\BBone_{h>a} \,
\Psi''(h) \, |\nabla h|^2$. Remark that we need $h>a$ for $1/h$ to be bounded (since $\Psi''(u)=2$
when $u\leq 2$) at least for $h<2$.

If $m \geq 1$ the situation is more delicate. But $$m \leq \int_{I^d} h d\bar{\mu} + a \leq
(1/\mu(I^d))+a$$ so that if we choose $R$ (the length of the edge of $I^d$) large enough we may
assume that $\mu(I^d)\geq 3/4$, choose $a=1/3$ so that $m \leq 5/3 < 2$. In other words on
$\{h>2\}$, $h/m \geq 6/5$. It follows that $\Psi(h)=\Psi(m \, \frac hm) \leq c \, \Psi(\frac hm)$
on $\{h>2\}$, for some constant $c$ (recall the form of $\Psi$). Furthermore $\BBone_{h>2} \leq
\BBone_{\frac hm > \frac 65}$ so that one more time we may apply \eqref{eqPsisob}, and conclude as
in the case $m \leq 1$.

We have thus shown the existence of some $C$ such that $$\int_{I^d} \, \Psi(h) \, \BBone_{h > 2}
\, d\mu \, \leq \, C \, \int_{\R^d} \, \Psi''(h) \, |\nabla h|^2 \, d\mu \, .$$ With the previous
result the proof is completed.
\end{proof}

\begin{remark}\label{remarqueentropydecay}
Since $\Psi(u)$ behaves like $u \log u$ at infinity, the previous result has the following
consequence : if $V$ has some exponential moment, then $$\int P^*_t h \, \log P_t^* h d\mu \, \leq
\, C \, e^{-\eta t} \, \left(1 \vee \int \Psi(h) \, \log_+(\Psi(h)) \, d\mu\right) \, \leq \, C'
\, e^{-\eta t} \, \left(1 \vee \int h \, \log_+^2(h) \, d\mu\right) \, .$$ This result is (at a
qualitative level) a little bit weaker than the one we obtain in this case in Theorem
\ref{thmMTlog}, since there we can replace the exponent 2 by any exponent greater than 1.

It should also be interesting to extend this kind of result to (strongly) hypoelliptic operators
as in Corollary \ref{corweight2}. The key would be to prove a local log-Sobolev inequality for the
corresponding $\Gamma$. We strongly suspect that some inequality of this type is true, but we did
not find any reference about it.
\end{remark}
\bigskip

%%%%%%%%%%%%%%%%%%%%
%%%
%%%%%%%%%%%%%%%%%%%%

\section{\bf Fully degenerate cases, towards hypocoercivity.}\label{sechypo}

Proposition \ref{proptrou} and Theorem \ref{thmentropydecaywithlyap} are hypocoercive results in
Villani's terminology. The former shows a coercivity property in $\L^2(W \mu)$ norm, which is
stronger than the $\L^2(\mu)$ norm, while the latter can be interpreted in terms of
semi-distances. We refer to \cite{Vilmad} for a nice presentation of hypocoercivity. In studying
fully degenerate cases, Villani introduces higher order functional inequalities (reminding the
celebrated $\Gamma_2$ criterion for logarithmic Sobolev inequality), see equation (11) in
\cite{Vilmad} and more generally \cite{Vilhypo}. These higher order inequalities enable him to
introduce Lie brackets of the diffusion vector fields with the drift vector field, hence are
clearly related to some hypoelliptic situation of H\"{o}rmander type. A deep study of the spectral
theory of hypoelliptic operators is done in \cite{HN}, and we refer to the references in both
\cite{HN,Vilhypo} for more details and contributors. Also notice that the hypocoercivity
phenomenon was first studied by H\'erau and Nier (see \cite{HN04}) by using pseudo-differential
calculus (also see some recent work by H\'erau on his Web page).
\medskip

Since the existence of a Lyapunov function does not immediately rely on non degeneracy, it is
natural to consider fully degenerate cases from this point of view. Note that Theorem
\ref{thmlyapgap} requires a local Poincar\'e inequality, hence is not adapted, while the method in
section \ref{MeynTweedie} furnishes some exponential decay for the variance but controlled by some
$\L^{p}$ norm.

In this section we shall recall the results in \cite{Vilhypo} for the particular example of the
kinetic Fokker-Planck equation. Then we shall see that this example enters the framework of
Meyn-Tweedie approach, following \cite{wudamp} and \cite{DFG} who indicated how to build some
Lyapunov function.

First we recall what the kinetic Fokker-Planck equation is. Let $F$ be a smooth function on
$\R^d$. We consider on $\R^{2d}$ the stochastic differential system ($x$ stands for position and
$v$ for velocity)
\begin{eqnarray}\label{eqdamping}
dx_t & = & v_t \, dt \\ dv_t & = & dB_t \, - \, v_t \, dt \, - \, \nabla F(x_t) \, dt \nonumber
\end{eqnarray}
associated with $$L = \frac 12 \Delta_v \, + \, v \, \nabla_x \, - \, (v \, + \, \nabla F(x)) \,
\nabla_v \, .$$ Define
\begin{equation}\label{eqmesinv}
\mu(dx,dv) = e^{- \, (|v|^2 + 2 F(x))} \, dx \, dv = e^{- H(x,v)} \, dx \, dv
\end{equation}
which is assumed to be a bounded measure (in the sequel we shall denote again by $\mu$ the
normalized (probability) measure $\mu/\mu(\R^{2d})$).

If $F$ is bounded from below, it is known that \eqref{eqdamping} has a pathwise unique, non
explosive solution starting from any $(x,v)$. Actually the statement in \cite{wudamp} Lemma 1.1 is
for a weak solution since Wu is using Girsanov theory. But introduce the stopping time $\tau_R=
\inf \{s\geq 0 ; |v_t|\geq R\}$. Since $|x_{t\wedge \tau_R}|\leq Rt +|x|$ pathwise uniqueness
holds up to each time $\tau_R$ and the explosion time is the limit of the $\tau_R$'s as $R$ goes
to infinity. That this limit is almost surely $+\infty$ is proved by Wu at the top of p.210 in
\cite{wudamp}.

Furthermore $\mu$ is in this case the unique invariant measure.

Let us make three additional remarks
\begin{itemize} \item $\mu$ is not symmetric, \item $L$ is fully degenerate, in particular since
$\Gamma f = |\nabla_v f|^2$ any function $f(x,v)=g(x)$ with $\int f d\mu =0$ is such that $\Gamma
f =0$ so that the Poincar\'e inequality (with $\Gamma$) is not true for $\mu$, \item the
Bakry-Emery curvature of the semi-group (see \cite{logsob} Definition 5.3.4) is equal to $-
\infty$.
\end{itemize}

The main results in \cite{Vilhypo} about convergence to equilibrium for this equation are
collected below
\begin{theorem}\label{thmvilkin}\textbf{Villani \cite{Vilhypo} Theorems 29,31,32}.
\begin{enumerate}
\item[(1)] \quad Define $H^1(\mu) := \{ f\in \L^2(\mu) ; \nabla f \in \L^2(\mu)\}$ equipped with
the semi-norm $\parallel f \parallel_{H^1(\mu)} = \parallel \nabla f\parallel_{\L^2(\mu)}$.

Assume that $|\nabla^2 F| \leq c (1+|\nabla F|)$ and that the marginal law $\mu_x(dx)=e^{-2F(x)}
dx$ satisfies the classical Poincar\'e inequality for all nice $g$ defined on $\R^d$
$$Var_{\mu_x}(g) \, \leq \, C \, \int_{\R^d} \, |\nabla g|^2(x) \, \mu_x(dx) \, .$$ Then there exist
$C$ and $\lambda$ positive such that for all $f \in H^1(\mu)$, $$\parallel P^*_t f \, - \, \int f
d\mu\parallel_{H^1(\mu)} \, \leq \, C \, e^{-\lambda t} \, \parallel f \parallel_{H^1(\mu)} \, .$$
\item[(2)] \quad With the same hypotheses, there exists $C$ such that for all $1\geq \varepsilon
>0$ and all $t>\varepsilon$, $$Var_\mu(P^*_t f) \, \leq
\, C \, \varepsilon^{-3/2} \, e^{-\lambda (t-\varepsilon)} \, Var_\mu(f) \, .$$ \item[(3)] \quad
Assume that $|\nabla^j F| \leq c_j$ for all $j\geq 2$ and that $\mu_x$ satisfies a (classical)
log-Sobolev inequality $$Ent_{\mu_x}(g^2) \, \leq \, C \, \int_{\R^d} \, |\nabla g|^2(x) \,
\mu_x(dx) \, .$$ Then for all $h\geq 0$ such that $\int h d\mu =1$ and satisfying
$$\forall k \in \N \quad , \quad \int \, (1+|x|+|v|)^k \, h(x,v)  \, d\mu \, < \, +\infty \, ,$$
it holds for some $\lambda > 0$, $$\int \, P^*_t h \, \log (P^*_t h) \, d\mu \, \leq \, C(h) \,
e^{- \lambda t}$$ where $C(h)$ depends on the above moments.
\end{enumerate}
\end{theorem}

It is worthwhile noticing that since $\mu$ is a product measure of $\mu_x$ and a gaussian measure,
$\mu$ inherits the classical Poincar\'e or log-Sobolev inequality as soon as $\mu_x$ satisfies one
or the other. Part (2) in the previous result is simply an hypoelliptic regularization property,
and some hypotheses can be slightly improved (see \cite{Vilhypo} Theorems 29,31 and 32 for the
details). However, it has to be noticed that $C>1$ (otherwise $\mu$ would satisfy a Poincar\'e
inequality with $\Gamma$) and that the Bakry-Emery curvature has to be $- \infty$ for the same
reason.
\medskip

In \cite{wudamp}, Wu gave some sufficient conditions for the existence of a Lyapunov
function for this (and actually more general) model (see \cite{wudamp} Theorem 4.1). We recall and
extend this result below. First define
\begin{equation}\label{eqwulyap}
\Lambda_{a,b}(x,v) = a H(x,v) + b(\langle v,\nabla G(x)\rangle + G(x))
\end{equation}
where $G$ is smooth, $a$ and $b$ being positive parameters.

\begin{theorem}\label{thmwuexp}
Assume that $F$ is bounded from below and that there exists some $G$ satisfying
\begin{enumerate}
\item[(1)] \quad $\liminf_{|x|\to +\infty} \langle \nabla G(x),\nabla F(x)\rangle = 2c >0$,
\item[(2)] \quad $\parallel \nabla^2 G\parallel_\infty \, < \, c/16d$ , \item[(3)] \quad there
exists $\kappa > 0$ such that for all $x$, $|\nabla G(x)|^2 \leq \kappa (1+|\langle\nabla
F(x),\nabla G(x)\rangle|)$, \item[(4)] \quad $\Lambda_{a,b}$ is bounded from below.
\end{enumerate}
Then for all $0<\varepsilon$ one can find a pair $(a,b)$ such that $\max(a,b)\leq \varepsilon$ for
which $V_{a,b}(x,v)= e^{\Lambda_{a,b}(x,v) - \inf_{x,v} \Lambda_{a,b}(x,v)}$ is a Lyapunov
function.

Hence if there exists $\eta>0$ such that $\int e^{\Lambda_{\eta,\eta}(x,v)} d\mu < +\infty$, for
each $p>1$ one can find a Lyapunov function $V_p \in \L^p(\mu)$, so that there exists $\lambda >0$
such that for each $q>2$ there exists $C_q$ such that
$$\Var_\mu(P_t^* f) \, \leq \, C_q \, e^{- (\frac{q-2}{q-1}) \, \lambda \, t} \, \parallel f -
\int f d\mu\parallel_q^2 \, .$$
\end{theorem}
\begin{proof}
 Elementary computation yields
$$LV_{a,b}/V_{a,b} = - 2a|v|^2(1-a)+ad+2ab\langle v,\nabla G\rangle+\frac 12 b^2 |\nabla G|^2 + b\langle
\nabla^2 G \, v,v\rangle - b\langle \nabla F, \nabla G \rangle \, .$$ Our aim is to choose $G$ for
the right hand side to be negative outside some compact set. A rough majorization gives
$$LV_{a,b}/V_{a,b} \leq  (-2a(1-a)+b |\nabla^2 G (x)|+4ab) |v|^2 - b \langle \nabla G,\nabla F\rangle +
(\frac{b^2}{2}+4ab) |\nabla G|^2 + ad \, .$$ We have thanks to (3) $$- b \langle \nabla G,\nabla
F\rangle + (\frac{b^2}{2}+4ab) |\nabla G|^2 + ad \leq b(-1 + \kappa (\frac b2 + 4a)) \langle
\nabla G,\nabla F\rangle + (ad+\kappa b(\frac b2 + 4a))$$ so that if we choose $a$ and $b$ small
enough for $\kappa (\frac b2 + 4a)\leq \frac 12$ the first term is less than $-cb$ for $|x|$ large
enough thanks to (1). Hence if we choose $ad+\kappa b(\frac b2 + 4a) < cb/2$ we get
$LV_{a,b}/V_{a,b} \leq -cb/2$ for $|x|$ large and all $v$ as soon as $$-2a(1-a)+b |\nabla^2 G
(x)|+4ab \leq 0 \, .$$ We may thus first choose $a$ and $b$ small enough for $\kappa (\frac b2 +
4a)<c/4$, so that it remains to choose $a<cb/4d$.

Now if $|x|\leq L$, $(LV_{a,b}/V_{a,b})(x,v) \to - \infty$ as $|v| \to +\infty$ as soon as
$$-2a(1-a)+b |\nabla^2 G (x)|+4ab < 0 \, .$$ We may choose $b\leq 1/8$ and $a\leq 1/2$ so that we
only have to check $-a/2 + b |\nabla^2 G (x)| < 0$, i.e. $a/2 > cb/16d$ thanks to (2). This is
possible since our unique constraint is $a/2<cb/8d$.

We have thus obtained the existence of a Lyapunov function for some pair $(a,b)$ with both $a$ and
$b$ as small as we want. This Lyapunov function thus belongs to $\L^p$ if $a$ and $b$ are small
enough, according to our integrability hypothesis. It remains to apply Theorem \ref{thmMT} to
conclude (all the other hypotheses in Theorem \ref{thmDFG} are satisfied here, see
\cite{wudamp,DFG} for the details).
\end{proof}

\begin{example}
Let us describe some examples.
\begin{enumerate}
\item[(1)] \, \, \, \textbf{(Wu \cite{wudamp})} \quad Assume the drift condition  $\liminf_{|x|\to
+\infty} \langle x,\nabla F(x)\rangle/|x| = 2c > 0$. Then we may choose $G(x)=|x|$ for $|x|$
large, and $|\nabla^2 G(x)| \leq \varepsilon$ for all $x$. This is the situation discussed in
\cite{wudamp}. Notice that $\mu_x$ satisfies a classical Poincar\'e inequality (see e.g. section
\ref{secex1}) so that the hypotheses of Theorem \ref{thmvilkin} are satisfied. \item[(2)] \quad A
little more general situation is for $F$ going to infinity, satisfying $$\liminf_{|x|\to +\infty}
|\nabla F(x)|^2 = 2c > 0 \textrm{ and }|\nabla^2 F|\ll |\nabla F| \textrm{ at infinity.}$$ In this
case also $\mu_x$ satisfies a classical Poincar\'e inequality as we saw in section \ref{secex1}
(if $d=1$ the converse is true). If $|\nabla^2 F(x)| \to 0$ as $|x| \to +\infty$ we may choose a
function $G$ such that $|\nabla^2 G(x)| \leq \varepsilon$ for all $x$ and $G(x)=F(x)$ for $x$
large. This function will satisfy all (1),(2),(3). For (4) and the integrability condition to be
satisfied it is enough to assume in addition that $$|\nabla F(x)|^2/F(x) \textrm{ goes to $0$ at
infinity.}$$ This is the case for $F(x)=|x|^p$ at infinity for $1\leq p < 2$. \item[(3)] \quad If
the latter condition is not satisfied we may take $G=F^\alpha$ for some $\alpha \leq 1$. But in
this situation we can obtain a better Lyapunov function and study convergence in entropy.
\end{enumerate}
\end{example}

\begin{remark}
The $\L^2$ convergence in Theorem \ref{thmvilkin} is optimal, hence we cannot expect to improve it
and actually the controls we obtained in Theorem \ref{thmwuexp} are weaker. In addition, in the
last version of his work (see \cite{vilgroshypo}) Villani gives some explicit bounds for the
constants involved. As we said, such estimates are not yet available in Theorem \ref{thmDFG}.

However, Villani's approach uses the classical Poincar\'e inequality in an essential way, and only
gives exponential decay results. Examples for the existence of $\phi$-Lyapunov functions for this
kinetic model are given in \cite{DFG} section 4.3 Indeed consider $F(x)\sim|x|^p$ for large $|x|$
with $0<p<1$.  Attentive calculations show that one can consider smooth $G$ with $\nabla G(x)=|x|^m$ for
large $|x|$ with $1-p<m\le1$,
$$e^{s\Lambda_{a,b}^\delta(x,v)- \inf_{x,v}s\Lambda_{a,b}^\delta(x,v)}\qquad (m+1)\delta\le p,$$
as a $\phi$-Lyapunov function for well chosen $s,a,b$, with $\phi(t)=t/\ln^{1/\delta-1}t$.
Combined with Theorem  \ref{thmMT} we thus get a subexponential decay in a situation where it is
known that there is no exponential decay, thanks to an argument by Wu \cite{wudamp}. We refer to
\cite{DFG} for the polynomial decay case. We shall not go further in this direction here, but
Theorems \ref{thmDFG} and \ref{thmMT} thus allow to study a larger field of potentials.
\end{remark}

As we said before we turn to the study of entropy decay.

This time we shall directly use $\Lambda_{a,b}+M=V_{a,b}$ as a Lyapunov function, for $M$ large
enough. Indeed
$$L V_{a,b}(x,v) = ad - 2a |v|^2 - b \langle \nabla F(x),\nabla G (x)\rangle + b \langle
\nabla^2 G(x) v,v\rangle \, .$$ Our aim is to find $G$ and $\eta > 0$ such that $LV_{a,b} \leq -
\eta V_{a,b}$ outside some compact set. We shall choose $G(x)=F^{1-\alpha}(x)$ for large $x$, for
some $0\leq \alpha<1$, assuming that $F$ is non-negative outside some compact set. Actually we
shall assume that $F$ goes to infinity at infinity. With all these choices
$$\Lambda_{a,b}(x,v) \geq a |v|^2 + 2a F(x) - b |v| \, \frac{|\nabla F(x)|}{F^{\alpha}(x)} $$ is
bounded from below as soon as $|\nabla F(x)|^2/F^{1+2\alpha}(x)$ goes to $0$ at infinity or if
this ratio is bounded and $b/a$ small enough.
\smallskip

Now if $\alpha > 0$, $$ \langle \nabla^2 G(x) v,v\rangle = (1-\alpha) F^{-\alpha}(x) \, \langle
\nabla^2 F(x) v,v\rangle - \alpha (1-\alpha) \, F^{-(1+\alpha)}(x) \, \langle \nabla F(x)
,v\rangle^2 \, ,$$ so that for $x$ large,
\begin{equation}\label{eqlyapentropy}
L V_{a,b}(x,v) \leq  ad - 2a |v|^2 - b (1-\alpha) F^{-\alpha}(x) \,|\nabla F(x)|^2 + b (1-\alpha)
F^{-\alpha}(x) \langle \nabla^2 F(x) v,v\rangle \, .
\end{equation}
To show that $V_{a,b}$ is a Lyapunov function, using the same majorization as in the proof of the
latter Theorem, it is enough to show that we can find some $\eta>0$ such that for $x$ large
\begin{equation}\label{eqlyapentropy2}
\left((2-\eta)a - 2 b \eta - b (1-\alpha)\frac{|\nabla^2 F(x)|}{F^{\alpha}(x)}\right) |v|^2 \, +
\, b \frac{|\nabla F(x)|^2}{F^{\alpha}(x)} \, \left(1 - \alpha - \frac{2
\eta}{F^{\alpha}(x)}\right) -
\end{equation}
$$ -
\left(M+ad + 2a \eta F(x)+b \eta F^{1-\alpha}(x)\right) \, \geq \, 0 \, .$$ Note that the same
result holds true for $\alpha = 0$.

The situation is now quite simple : first we shall assume that $|\nabla F(x)|^2 \geq \kappa
F^{1+\alpha}(x)$ for large $x$, so that for any $b$ we may choose $\eta$ small enough for the sum
of the last two terms to be positive; next we have to assume that $|\nabla^2 F(x)|/F^{\alpha}(x)$
is bounded, so that we may choose $b$ small enough for the coefficient of $|v|^2$ to be positive.
Of course for $|x|\leq L$ \eqref{eqlyapentropy2} has to be replaced by the correct one involving
$G$, but $G$ being smooth it is enough again to choose $b$ and $\eta$ small enough.

Choosing $a$ small enough we see that $\int e^{pV_{a,b}} d\mu < +\infty$, so that applying Theorem
\ref{thmMTlog} and H\"{o}lder-Orlicz inequality to bound $\int Vh d\mu$ we have obtained
\begin{theorem}\label{thmentropyBCG}
Assume that $F(x) \to +\infty$ as $|x| \to +\infty$ (hence is bounded from below) and that there
exists $0\leq \alpha < 1$ such that the following holds
\begin{enumerate}
\item[(1)] \quad there exist $c$ and $C$ such that for $|x|$ large, $$c \, F^{1+\alpha}(x) \leq
|\nabla F(x)|^2 \, \leq C \, F^{1+2\alpha}(x) \, ,$$ \item[(2)] \quad $|\nabla^2
F(x)|/F^{\alpha}(x)$ is bounded (for $|x|$ large).
\end{enumerate}
Then for all $p>1$ one can find a Lyapunov function $V_p$ such that $\int e^{pV} d\mu < +\infty$.
Hence there exists $\lambda>0$ such that for any $1>\beta>0$ there exists $C_\beta$ such that for
all density of probability $h$,
$$\int \, P_t^* h \, \log P^*_t h \, d\mu \leq \, C_\beta \, e^{- \beta \, \lambda \, t} \,
\left(1+ \int h \log h d\mu \right)^\beta \, \left(\int |h-1| |\log h|^{\frac{1}{1-\beta}}
d\mu\right)^{1-\beta} \, .$$
\end{theorem}

\begin{example}
If $F(x)=|x|^p$ for some $p\geq 2$ and large $|x|$, then we may apply the previous Theorem with
$$\frac{p-2}{2p} \, \leq \, \alpha \, \leq \, \frac{p-2}{p} \, .$$
\end{example}

\begin{remark}\label{remimpo}
As it is shown in \cite{Cat5} the condition $|\nabla F(x)|^2 \geq \eta F(x) + \Delta F(x)$ for
large $x$ implies a classical logarithmic Sobolev inequality for $\mu$. Hence if $|\nabla^2 F|
\leq C(1+\nabla|F|)$ our hypothesis (1) in Theorem \ref{thmentropyBCG} implies a classical
logarithmic Sobolev inequality, as it is asked in Theorem \ref{thmvilkin} (3).

But case (3) in Theorem \ref{thmvilkin} is (very) roughly the case where $c |x|^2 \leq F(x) \leq C
\, |x|^2$ for some positive $c$. Our result covers more ``convex at infinity'' cases.

Finally, even if we do not have explicit constants, our hypotheses on $h$ seem to be weaker than
the moment conditions in Theorem \ref{thmvilkin}. For instance if $F(x)=|x|^2/2$ we may choose with $a>0$
$$h(x,v)= \frac{e^{|x|^2+|v|^2}}{\left(1+|x|^{d+1}+|v|^{d+1}\right)^{a+1}}$$ for any $\beta<1-2/(a(d+1))$,
while this $h$ does not
fulfill the hypotheses of Theorem \ref{thmvilkin} (3) (requires all $\beta<1$!).
\end{remark}

\begin{remark}
Of course, since for any density of probability $h$ it holds $\int h \log h d\mu := Ent_\mu(h)
\leq Var_\mu(h)$, the relative entropy is decaying at least with the same rate as the variance,
hence Theorem \ref{thmwuexp} furnishes some decay. The study of relative entropy in \cite{HN04} is
based on this argument.
\end{remark}

\begin{remark}
Remark that the generator $L$ can be written in H\"{o}rmander's form $L=\frac 12 \, X_1^2 + X_0$
where the vector fields $X_i(x,v)$ are given by $X_1(x,v)= \partial_v$ and $X_0(x,v)= v
\partial_x - (v+\nabla F(x)) \partial_v$. Hence the Lie bracket $[X_1,X_0](x,v)=\partial_x -
\partial_v$ is such that $X_1$ and $[X_1,X_0]$ generate the tangent space at any $(x,v)$.
Furthermore $|X_1|^2+|[X_1,X_0]|^2$ is uniformly bounded from below by a positive constant. Hence
Malliavin calculus allows us to show that, for any $t>0$, the law of $(x_t,v_t)$ starting from any
point $(x,v)$ has a $C^{\infty}$ density $p_t$ w.r.t Lebesgue measure, hence a smooth density
$h_t$ w.r.t. $\mu$. Furthermore $p_t$ satisfies some gaussian upper bound. However we do not know
how to show that $h_t \in \L^2(\mu)$. The latter is shown in \cite{HN04}, but starting with some
particular initial absolutely continuous laws. Due to the gaussian part of $\mu$, exponent 2 is
optimal for such a result.
\end{remark}

\bigskip

\bigskip
\bibliographystyle{plain}

\end{document}